\documentclass[moor]{informs1}              

\usepackage{amsmath, amsfonts, amssymb,  mathrsfs, stmaryrd, comment}
\newcommand{\ab}[1]{\langle #1\rangle}

\newcommand{\lbar}[1]{\underline{#1}}
\newcommand{\ubar}[1]{\overline{#1}}

\newcommand{\norm}[1]{\| #1\|}

\usepackage{natbib}
 \NatBibNumeric
 \def\bibfont{\small}%
 \bibpunct[, ]{[}{]}{,}{n}{}{,}%

\usepackage[colorlinks=true,breaklinks=true,bookmarks=true,urlcolor=blue,
     citecolor=blue,linkcolor=blue,bookmarksopen=false,draft=false]{hyperref}

\def\EMAIL#1{\href{mailto:#1}{#1}}

\TheoremsNumberedThrough     

\EquationsNumberedThrough    

\begin{document}


\RUNAUTHOR{Bayraktar and Zhang}

\RUNTITLE{FTAP under transaction costs and model uncertainty}

\TITLE{Fundamental Theorem of Asset Pricing under Transaction costs and Model uncertainty}

\ARTICLEAUTHORS{
\AUTHOR{Erhan Bayraktar}
\AFF{Department of Mathematics, University of Michigan, 530 Church Street, Ann Arbor, MI 48104, \EMAIL{erhan@umich.edu}}
\AUTHOR{Yuchong Zhang}
\AFF{Department of Statistics, Columbia University, 1255 Amsterdam Avenue, New York, NY 10027, \EMAIL{yz2915@columbia.edu}}
}

\ABSTRACT{
We prove the Fundamental Theorem of Asset Pricing for a discrete time financial market where trading is subject to proportional transaction cost and the asset price dynamic is modeled by a family of probability measures, possibly non-dominated.
Using a backward-forward scheme, we show that when the market consists of a money market account and a single stock, no-arbitrage in a quasi-sure sense is equivalent to the existence of a suitable family of consistent price systems. We also show that when the market consists of multiple dynamically traded assets and satisfies \emph{efficient friction}, strict no-arbitrage in a quasi-sure sense is equivalent to the existence of a suitable family of strictly consistent price systems.
}


\MSCCLASS{Primary: 60G42, 91B28; secondary: 03E15}
\ORMSCLASS{Primary: finance: asset pricing; secondary: probability}
\HISTORY{Received on May 13, 2014; Second version received on May 8, 2015; Accepted on August 29, 2015.}

\maketitle

%

\section{Introduction}
The Fundamental Theorem of Asset Pricing (FTAP), as suggested by its name, is one of the most important theorems in mathematical finance and has been established in many different settings: discrete and continuous, with and without transaction costs. It relates no-arbitrage concepts to the existence of certain fair pricing mechanisms, which provides the rationale for why in duality theory, it is often reasonable to assume that the dual domain is non-empty. 

We consider a discrete time, finite horizon financial market consisting of stocks and a zero-interest money market account. When the market is frictionless and modeled by a single probability measure, the classical result by Dalang-Morton-Willinger \cite{DMW90} asserts there is no-arbitrage if and only if the set of equivalent martingale measures is not empty. With proportional transaction costs, the set of martingale measures is replaced by the set of consistent price systems (CPS's) or strictly consistent price systems (SCPS's). Equivalence between no-arbitrage and existence of a CPS is established by Kabanov and Stricker \cite{KS01} for finite probability space $\Omega$, and by Grigoriev \cite{Grigoriev05} when the dimension is two. Such equivalence in general does no hold for infinite spaces and higher dimensions (see Section 3 of Schachermayer \cite{Schachermayer04} and page 128-129 of Kabanov and Safarian \cite{Tran.Cost} for counter examples). For such an equivalence one needs the notion of robust no-arbitrage introduced by Schachermayer in \cite{Schachermayer04}, where he showed that this notion is equivalent to the existence of an SCPS. Alternatively, robust no-arbitrage can be replaced by strict no-arbitrage plus efficient friction (see Kabanov et al. \cite{KRS02}). There exist a few different proofs of the FTAP under transaction costs. Besides the proofs of Schachermayer \cite{Schachermayer04} and Kabanov et al. \cite{KRS02} which rely on the closedness of the set of hedgeable claims and a separation argument, there is a utility-based proof by Smaga \cite{Smaga12} and proofs based on random sets by Rokhlin \cite{Rokhlin07b}. 

In recent years, model uncertainty has gained a lot of interest in financial mathematics since it corresponds to a more realistic modeling of the financial market. By model uncertainty, we mean a family $\mathcal{P}$ of probability measures, usually non-dominated. Each member of $\mathcal{P}$ represents a possible model for the stock. One should think of $\mathcal{P}$ as being obtained from market data. We have a collection of measures rather than a single one because we do not have point estimates but confidence intervals. The non-dominated case is generally much harder because the classical separation argument does not work. In a market model without transaction cost, the recent work by Bouchard and Nutz \cite{BN13} used a local analysis to establish the equivalence between the absence of arbitrage in a quasi-sure sense and the existence of an ``equivalent" family of martingale measures. Acciaio et al. \cite{Schachermayer13} obtained a different version of the FTAP by working with a different no-arbitrage condition which excludes model-independent arbitrage, i.e. arbitrage in a sure sense only. 

In this paper, we prove the FTAP when both proportional transaction cost and model uncertainty are present. We start with a financial market consisting of a money market account and a single stock and introduce the no-arbitrage concept NA$(\mathcal{P})$. We prove its equivalence to the existence of a family of CPS's, and later extend our results to a market with multiple assets where NA$(\mathcal{P})$ is replaced by NA$^s(\mathcal{P})$ plus efficient friction, and CPS is replaced by SCPS. Similar to the classical theory (in the set-up of which there is a dominating measure), the two-dimensional case is different from the multi-dimensional case in that no-arbitrage alone is sufficient to imply the existence of a CPS. That is the market model need not satisfy the efficient friction hypothesis. Hence the two-dimensional case is worth a separate study, and will be the main focus of this paper. On a related note, when there is a single stock, Dolinsky and Soner \cite{DS14} proved the super-hedging theorem (by first discretizing the state space and then taking a limit) and stated the FTAP as a corollary. We follow a different methodology and state the result for an arbitrary collection of probability measures instead of the collection of all probability measures as they do. Thanks to our method we are able to work with a more general structure on the proportional transaction cost instead of taking it a constant, which is useful as this proportion in real markets is likely to change with changing market conditions over time. Our contribution can be seen as an extension of Bouchard and Nutz \cite{BN13} to the transaction cost case, and a generalization of the results of Rokhlin \cite{Rokhlin07b} on the martingale selection problem to the non-dominated case.

Similar to Bouchard and Nutz \cite{BN13}, we proceed in a local fashion: first obtain a CPS for each single-period model and then do pasting using a suitable measurable selection theorem. Although our proof follows the ideas in \cite{BN13}, the multi-period case turns out to be quite different when transaction costs are added. A distinct feature for frictionless markets is that the absence of arbitrage for the multi-period market is equivalent to the absence of arbitrage in all single-period markets. So it is enough to look at each single period separately and paste the martingale measures together. This equivalence, however, breaks down in the presence of transaction cost. A simple example is the following two-period market: $\lbar{S}_0=1, \ubar{S}_0=3, \lbar{S}_1=2, \ubar{S}_1=4, \lbar{S}_2=3.5, \ubar{S}_2=5$ where an underline denotes the bid price and an overline denotes the ask price. Each period is arbitrage-free, but buying at time 0 and selling at time 2 is an arbitrage for the two-period market. So we cannot in general paste two one-period martingale measures to get a two-period martingale measure; in particular, the endpoints of the underlying martingales constructed for each single period may not match. We need to solve a non-dominated martingale selection problem. The martingale selection problem when $\mathcal{P}$ is a singleton was studied by Rokhlin in a series of papers \cite{Rokhlin05, Rokhlin07a, Rokhlin07b} using the notion of support of regular conditional upper distribution of set-valued maps. In our case, it is difficult to talk about conditional distribution due to the lack of dominating measure. Nevertheless, we got some inspiration from Rokhlin \cite{Rokhlin07a,Rokhlin07b} and Smaga \cite{Smaga12} and developed a backward-forward scheme:

\emph{Backward recursion}: modify the original bid-ask prices backward in time by potentially more favorable ones to account for the missing future investment opportunities; 

\emph{Forward extension}: extend the CPS forward in time in the modified market.


When there is no dominating measure, the backward-forward scheme brings some measurability issues. On one hand, neither Borel nor universal measurability is preserved under the backward recursion. This can be circumvented by working with lower semi-analytic bid prices and upper semi-analytic ask prices. On the other hand, the proofs of Lemmas 4.6 and 4.8 in Bouchard and Nutz \cite{BN13} rely critically on the stock prices being Borel measurable. However, our modified market only has semi-analytic bid-ask prices, which makes measurable selection in many places challenging (see Remarks \ref{rmk:Nt} and \ref{rmk:pair}). To overcome this difficulty, we break up the selection of a CPS into two steps: first select a pair of measures $(Q_1, Q_2)$ using the semi-analyticity of the bid and ask prices separately, and then find a convex weight function $\lambda$ which plays an important role in boundary extension, which is an additional step that needs to be taken care of in dimension two because of the lack of the efficient friction hypothesis.  A measurable one-period CPS can be constructed out of the measurable triplet $(Q_1, Q_2, \lambda)$ and the bid-ask prices of the modified market.

Similar measurability issues exist in the multi-dimensional case. Fortunately, it turns out that doing backward recursion on the dual cones of the solvency cones preserves graph analyticity which is enough to make measurable selection in the forward extension possible. Here at each step in the forward extension, we do not select a one-period SCPS $(Q, Z)$ with a prescribed initial value $Z_0=z$ directly, but a vector of equivalent measures $(Q,\mu^1, \ldots, \mu^d)$ related to $(Q, Z)$ via $Z^i_{1}=z^i d\mu^i/dQ$. This allows us to use the fact that there exists a jointly Borel measurable function $(Q,\mu^i,\omega')\mapsto \frac{d\mu^i}{dQ}(\omega')$, which is key to measurable selection. 

We point out that during the review process, a paper by Bouchard and Nutz \cite{BN14} which is closely related to ours came out. There the authors propose a different notion of no-arbitrage (a quasi-sure version of no-arbitrage of the second kind), denoted by NA$_2(\mathcal{P})$, with which they proved a multi-dimensional version of the FTAP under the extra assumption of efficient and bounded friction. We find their no-arbitrage notion a bit too strong. It means that the market is already in a good form for martingale extension, so that the backward recursion is avoided. A simple one-period, single stock market $\lbar{S}_0=1, \ubar{S}_0=3, \lbar{S}_1=2, \ubar{S}_1=4$ which satisfies NA$^s(\mathcal{P})$, thus reasonable in our opinion, fails NA$_2(\mathcal{P})$. Despite that NA$_2(\mathcal{P})$ fails for the original market, it holds for the corresponding modified market. The multi-dimensional part of our paper is a generalization of \cite{BN14} to allow for more general market structure. Apart from the no-arbitrage concept used, our paper also differs in two other ways. First, we do not assume bounded friction and $K^\ast_t\cap \partial\mathbb{R}^d_+=\{0\}$; efficient friction can also be dropped in the two-dimensional case. Second, our measurability assumption is weaker; we are able to handle markets with universally measurable solvency cones as long as their dual cones have analytic graphs.

The rest of the paper is organized as follows. Sections~\ref{sec:set-up}, \ref{sec:one_period} and \ref{sec:multi} are devoted to the two-dimensional case. In Section~\ref{sec:set-up}, we introduce the probabilistic framework, set up the financial market model, and state the main result. In Section~\ref{sec:one_period}, we prove some one-period results which serve as the building blocks. In Section \ref{sec:multi}, we present the backward-forward schemes for a multi-period market and prove the main theorem. Section~\ref{sec:multi-d} extends our methodology and results to the multi-dimensional case. A few technical lemmas are collected in Appendix.

\section{The Financial Market Model and Main Results}\label{sec:set-up}
We follow the notation in Bouchard and Nutz \cite{BN13}. Let $\mathbf{T}\in\mathbb{N}$ be the time horizon. Let $\Omega_1$ be a Polish space and $\Omega_t:=\Omega_1^t$ be the $t$-fold Cartesian product with the convention that $\Omega_0$ is a singleton. Denote by $\mathcal{B}(\Omega_t)$ the Borel sigma-algebra on $\Omega_t$, and by $\mathcal{F}_t$ the universal completion of $\mathcal{B}(\Omega_t)$. We write $(\Omega,\mathcal{F})$ for $(\Omega_\mathbf{T},\mathcal{F}_\mathbf{T})$. When we say a stochastic process is predictable or adapted, we mean with respect to the filtration $(\mathcal{F}_t)$. Let $\mathfrak{P}(\Omega_1)$ denote the set of all probability measures on $(\Omega_1, \mathcal{B}(\Omega_1))$. For each $t\in\{0,\ldots, \mathbf{T}-1\}$ and $\omega\in\Omega_t$, we are given a nonempty convex set $\mathcal{P}_t(\omega)\subseteq\mathfrak{P}(\Omega_1)$, representing the set of possible models for the $(t+1)$-th period. We assume $\mathcal{P}_t$ (considered as set-valued maps from $\Omega_t$ to $\mathfrak{P}(\Omega_1)$) has an analytic graph. This assumption ensures that $\mathcal{P}_t$ admits a universally measurable selector. Define the uncertainty set $\mathcal{P}\subseteq\mathfrak{P}(\Omega)$ of the multi-period market by
\begin{equation*}
\mathcal{P}:=\{P_0\otimes\cdots\otimes P_{\mathbf{T}-1}: \text{ each } P_t \text{ is a universally measureable selector of }\mathcal{P}_t\},
\end{equation*}
where 
\[P_0\otimes\cdots\otimes P_{\mathbf{T}-1}(A)=\int_{\Omega_1}\cdots\int_{\Omega_1} 1_A(\omega_1,\ldots,\omega_\mathbf{T})P_{\mathbf{T}-1}(\omega_1,\ldots,\omega_{\mathbf{T}-1};d\omega_\mathbf{T})\cdots P_0(d\omega_1), \ A\in\mathcal{F}.\]

Consider a financial market consisting of a money market account with zero interest rate, and a stock with strictly positive bid price $\lbar{S}_t$ and ask price $\ubar{S}_t$. $\lbar{S}_t, \ubar{S}_t:\Omega_t\rightarrow\mathbb{R}$ are assumed to be lower semi-analytic (l.s.a.) and upper semi-analytic (u.s.a.), respectively, for all $t$. A \textit{self-financing}\footnote{We allow agents to throw away non-negative quantities of the assets.} portfolio process is an $\mathbb{R}^2$-valued predictable process $\phi=(\phi^0,\phi^1)$ satisfying $\phi_0=0$ and $\triangle \phi_{t+1}^0\leq -(\triangle \phi_{t+1}^1)^+ \ubar{S}_t+(\triangle \phi_{t+1}^1)^-\lbar{S}_t$ for all $t=0,\ldots, \mathbf{T}$, where $\triangle \phi_{t+1}:=\phi_{t+1}-\phi_t$. Denote by $\mathcal{H}$ the set of self-financing portfolio processes. Let $A_\mathbf{T}:=\{\phi_{\mathbf{T}+1}:\phi\in\mathcal{H}\}$.  $A_\mathbf{T}$ is interpreted as the set of hedgeable claims (in terms of physical units) from zero initial endowment.

\begin{definition}\label{defNA}
We say NA($\mathcal{P}$) holds if for all $f\in A_\mathbf{T}$, $f\geq 0$ $\mathcal{P}$-quasi-surely (q.s.)\footnote{A set is $\mathcal{P}$-polar if it is $P$-null for all $P\in\mathcal{P}$. A property is said to hold $\mathcal{P}$-q.s. if it holds outside a $\mathcal{P}$-polar set.}
 implies $f=0$ $\mathcal{P}$-q.s..
\end{definition}

\begin{remark}
NA($P$) $\forall P\in\mathcal{P}$ implies NA($\mathcal{P}$). Indeed, let $f\in A_\mathbf{T}$ be such that $f\geq 0$ $\mathcal{P}$-q.s. hence $P$-a.s. for all $P\in\mathcal{P}$. For each $P$, NA($P$) implies $f=0$ $P$-a.s.. Since this holds for all $P\in\mathcal{P}$, $f=0$ $\mathcal{P}$-q.s. and NA($\mathcal{P}$) holds. The reverse direction is not true. Consider a one-period market with $S_0=2$ and $S_1(\omega)=1$, $S_1(\omega')=3$. Let $P_1=\delta_\omega, P_2=\delta_{\omega'}$ be the Dirac measures concentrated on $\omega$ and $\omega'$, respectively. Then it is easy to see that there is arbitrage under both $P_1$ and $P_2$, but not under $\mathcal{P}:=conv\{P_1,P_2\}$.\footnote{``conv" stands for convex hull.}
\end{remark}

\begin{definition}
A pair $(Q,\tilde{S})$ is called a consistent price system (CPS) if $\tilde{S}$ is a $Q$-martingale in the filtration $(\mathcal{F}_t)$ and $\tilde{S}_t\in[\lbar{S}_t,\ubar{S}_t]$ $Q$-a.s.. Denote the set of all consistent price systems by $\mathcal{Z}$.
\end{definition}

The main theorem of this paper is given below.
\begin{theorem}\label{FTAP}
The following are equivalent:
\begin{itemize}
\item[(i)] NA$(\mathcal{P})$ holds.
\item[(ii)] $\forall P\in\mathcal{P}, \exists (Q,\tilde{S})\in\mathcal{Z} \text{ such that }P\ll Q\lll\mathcal{P}$.\footnote{The notation $Q\lll\mathcal{P}$ is taken from \cite{BN13}. It means $Q\ll P$ for some $P\in\mathcal{P}$.}
\end{itemize}
\end{theorem}

\begin{remark}
Let $\mathcal{Q}$ be the collection of the first components of $\mathcal{Z}$ that are dominated by some measure in $\mathcal{P}$. Theorem \ref{FTAP}(ii) is equivalent to saying $\mathcal{P}, \mathcal{Q}$ are equivalent in terms of polar sets. When $\mathcal{P}$ is a singleton, we recover the classical result of the existence of an equivalent measure.
\end{remark}

\section{The One-Period Case}\label{sec:one_period}
In this section, we prove the more difficult direction of the FTAP (i.e. no-arbitrage implies the existence of CPS's) for a one-period market. To prepare for multi-period case, we also discuss how to construct martingales with prescribed initial values.

Let $(\Omega,\mathcal{F})$ be a measurable space with filtration $(\mathcal{F}_0,\mathcal{F}_1)$ and $\mathcal{F}_0=\{\emptyset,\Omega\}$. Let $\mathcal{P}\subseteq \mathfrak{P}(\Omega)$ be a nonempty convex set. The bid and ask price processes of the stock are given by constants $\lbar{S}_0,\ubar{S}_0$ and $\mathcal{F}_1$-measurable random variables $\lbar{S}_1, \ubar{S}_1$, respectively. 
Note that NA$(\mathcal{P})$ for this one-period market can be stated in the equivalent form: $\forall y\in\mathbb{R}$, $y^+(\lbar{S}_1-\ubar{S}_0)-y^-(\ubar{S}_1-\lbar{S}_0)\geq 0$ $\mathcal{P}$-q.s. implies $y^+(\lbar{S}_1-\ubar{S}_0)-y^-(\ubar{S}_1-\lbar{S}_0)= 0$ $\mathcal{P}$-q.s..

For each $P\in\mathcal{P}$, define
\begin{equation*}
\Theta_P:=\{R\in\mathfrak{P}(\Omega): P\ll R\lll\mathcal{P}, E^R[|\ubar{S}_1-\lbar{S}_0|+|\lbar{S}_1-\ubar{S}_0|]<\infty \}.
\end{equation*}
$\Theta_P$ is nonempty by Lemma \ref{Aequivmeas}.

\begin{lemma}\label{fundamental_lemma}
Suppose $\exists P_1,P_2\in\mathcal{P}$ satisfying $P_1(\ubar{S}_1-\lbar{S}_0>0)>0$ and $P_2(\lbar{S}_1-\ubar{S}_0<0)>0$. Then $\forall P\in\mathcal{P}$, $\exists Q_1, Q_2\in \Theta_P$ such that $Q_1\sim Q_2$ and
\begin{equation*}
E^{Q_1}[\ubar{S}_1-\lbar{S}_0]>0, \quad E^{Q_2}[\lbar{S}_1-\ubar{S}_0]<0.
\end{equation*}
\end{lemma}
\proof{Proof.}
Let $P\in\mathcal{P}$. We first show, along the same lines in the first paragraph on page 13 of \cite{BN13}, that there exist $P'_1, P'_2\in\Theta_P$ (not necessarily equivalent) such that 
\[E^{P'_1}[\ubar{S}_1-\lbar{S}_0]>0, \quad E^{P'_2}[\lbar{S}_1-\ubar{S}_0]<0.\]
Let $A:=\{\ubar{S}_1>\lbar{S}_0\}$.  We have $P_1(A)>0$. Define $R_1:=(P_1+P)/2$, use Lemma \ref{Aequivmeas} to replace $R_1$ by $R_2\sim R_1$ such that $R_2\in\Theta_P$, and further replace $R_2$ by $P'_1\sim R_2$ defined by $dP'_1/dR_2:=(1_A+\varepsilon)/E^{R_2}[1_A+\varepsilon]$. It can be checked that $P'_1\in\Theta_P$ and $E^{P'_1}[\ubar{S}_1-\lbar{S}_0]>0$ for $\varepsilon$ small enough. Similar, we obtain the existence of $P'_2$.

It remains to replace $P'_1$, $P'_2$ by $Q_1$, $Q_2$ with the additional requirement that $Q_1\sim Q_2$. To this end, let $Q_\lambda:=\lambda P'_1+(1-\lambda) P'_2$. It is easy to see that $\{Q_\lambda: \lambda\in (0,1)\}$ is a set of equivalent measures contained in $\Theta_P$. Moreover,  $E^{Q_\lambda}[\ubar{S}_1-\lbar{S}_0]>0$ for $\lambda$ sufficiently close to 1, and $E^{Q_\lambda}[\lbar{S}_1-\ubar{S}_0]<0$ for $\lambda$ sufficiently close to zero.\Halmos
\endproof

\begin{proposition}\label{FTAP1periodNA}
Suppose NA$(\mathcal{P})$ holds. Then $\forall P\in\mathcal{P}, \exists (Q,\tilde{S})\in\mathcal{Z} \text{ such that }P\ll Q\lll\mathcal{P}$.
\end{proposition}
\proof{Proof.}
Let $P\in\mathcal{P}$ and consider three cases:

Case 1. $\ubar{S}_1-\lbar{S}_0\leq 0$ $\mathcal{P}$-q.s.. In this case, NA$(\mathcal{P})$ implies $\ubar{S}_1-\lbar{S}_0=0$ $\mathcal{P}$-q.s.. We choose the CPS to be $Q:=P$, $\tilde{S}_0:=\lbar{S}_0$ and $\tilde{S}_1:=\ubar{S}_1$.

Case 2. $\lbar{S}_1-\ubar{S}_0\geq 0$ $\mathcal{P}$-q.s.. In this case, NA$(\mathcal{P})$ implies $\lbar{S}_1-\ubar{S}_0=0$ $\mathcal{P}$-q.s.. We choose the CPS to be $Q:=P$, $\tilde{S}_0:=\ubar{S}_0$ and $\tilde{S}_1:=\lbar{S}_1$.

Case 3. $\exists P_1,P_2\in\mathcal{P}$ such that $P_1(\ubar{S}_1-\lbar{S}_0>0)>0$ and $P_2(\lbar{S}_1-\ubar{S}_0<0)>0$. In this case, let $Q_1, Q_2\in\Theta_P$ be given by Lemma \ref{fundamental_lemma}. We can find $\lambda\in(0,1)$ such that 
\[\lambda E^{Q_1}[\ubar{S}_1-\lbar{S}_0]+(1-\lambda)E^{Q_2}[\lbar{S}_1-\ubar{S}_0]=0.\]
Define 
\begin{equation}\label{eq:Q}
Q:=\lambda Q_1+(1-\lambda) Q_2
\end{equation} 
and
\begin{equation}\label{eq:Stilde}
\tilde{S}_0:=(1-\lambda)\ubar{S}_0+\lambda\lbar{S}_0, \quad \tilde{S}_1:=\lambda \frac{dQ_1}{dQ}\ubar{S}_1+(1-\lambda)\frac{dQ_2}{dQ}\lbar{S}_1.
\end{equation}
We have $Q\in\Theta_P$, $\tilde{S}\in [\lbar{S}_1, \ubar{S}_1]$ $Q$-a.s., and $E^Q[\tilde{S}_1-\tilde{S}_0]=\lambda E^{Q_1}[\ubar{S}_1-\lbar{S}_0]+(1-\lambda)E^{Q_2}[\lbar{S}_1-\ubar{S}_0]=0$.\Halmos\endproof

\begin{remark}
We can directly use $P'_1, P'_2$ (defined in the proof of Lemma  \ref{fundamental_lemma}) instead of $Q_1, Q_2$ to construct the CPS in Case 3. The equivalence of $Q_1$ and $Q_2$ only matters in the multi-period case where it is important for us to construct martingales that stay away from the boundary of the bid-ask spread whenever possible (so that they are extendable to the next period). If $Q_1\sim Q_2$, $\tilde{S}$ defined by \eqref{eq:Stilde} will satisfy $\tilde{S}_1\in ri[\lbar{S}_1, \ubar{S}_1]$ $Q$-a.s..\footnote{``ri" stands for relative interior.}
\end{remark}

Recall that when we go to the multi-period case, we cannot directly paste two single-period CPSs, but to first make sure the starting point of the current-period martingale matches the terminal point of its parent. In other  words, we are interested in constructing martingales with certain prescribed initial values. Proposition \ref{1pmartingale} gives a set of starting points that admit a martingale extension.

For a random variable $S:\Omega\rightarrow \mathbb{R}$ and a nonempty family $\mathcal{R}$ of probability measures on $\Omega$, the support of the distribution of $S$ under $\mathcal{R}$, denoted by supp$_{\mathcal{R}}S$, is the smallest closed set $A\subseteq\mathbb{R}$ such that $P(S\in A)=1 \ \forall P\in\mathcal{R}$. Equivalently, a point $y$ belongs to $\text{supp}_{\mathcal{R}}S$ if and only if every open ball around $y$ has positive measure under some member of $\{P\circ S^{-1}: P\in\mathcal{R}\}$.

\begin{proposition}\label{1pmartingale} 
Suppose $\inf\text{supp}_\mathcal{P}\lbar{S}_1<s<\sup\text{supp}_\mathcal{P}\ubar{S}_1$. Then $\forall P\in\mathcal{P}$, $\exists Q_1, Q_2\in \Theta_P$ such that $Q_1\sim Q_2$ and
\[E^{Q_2}[\lbar{S}_1]<s<E^{Q_1}[\ubar{S}_1].\]
\end{proposition}
\proof{Proof.}
We can find $x\in\text{supp}_{\mathcal{P}}\lbar{S}_1$ and $y\in\text{supp}_{\mathcal{P}}\ubar{S}_1$ such that $x<s<y$. By definition of support, $\exists P_1, P_2\in\mathcal{P}$ satisfying $P_1(\ubar{S}_1>s)>0$ and $P_2(\lbar{S}_1<s)>0$. Lemma \ref{fundamental_lemma} applied to the market $\{s,[\lbar{S}_1, \ubar{S}_1]\}$ yields the desired $Q_1$ and $Q_2$.\Halmos
\endproof

\begin{remark}\label{rmk:1p-pair}
It is possible to show that $\forall s\in(\inf\text{supp}_\mathcal{P}\lbar{S}_1,\sup\text{supp}_\mathcal{P}\ubar{S}_1)$ and $P\in\mathcal{P}$, $\exists Q\in\Theta_P$ and $\lambda\in(0,1)$ such that $E^Q[(1-\lambda)\lbar{S}_1+\lambda\ubar{S}_1]=s$. For the proof, simply consider three cases: 1) $\exists R\in{\Theta_P}$ with $E^R[\lbar{S}_1]>s$, 2) $\exists R\in{\Theta_P}$ with $E^R[\ubar{S}_1]<s$, and 3) $\forall R\in{\Theta_P}$, $E^R[\lbar{S}_1]\leq s\leq E^R[\ubar{S}_1]$. 
We state the one-period result in terms of a pair of measures $(Q_1, Q_2)$ because it is more amenable for measurable selection when moving to the multi-period case (see Remark \ref{rmk:pair}). Once we have a pair of measures, a CPS can be easily constructed by means of \eqref{eq:Q} and \eqref{eq:Stilde}.
\end{remark}

\section{The Multi-period Case}\label{sec:multi}
In this section, we prove the FTAP for a multi-period market through a backward-forward scheme. Back to the setup in introduction, the set $\mathcal{P}$ is defined as the product of the nonempty convex sets $\mathcal{P}_t(\cdot)$ which have analytic graphs. Throughout this section, we also assume $\lbar{S}_t(\cdot)$ is l.s.a. and $\ubar{S}_t(\cdot)$ is u.s.a.. The reason for working with semi-analytic price processes is that we need a property that can be preserved under the backward recursion \eqref{XY_recursion}. Both Borel and universal measurability are not preserved (see Remark \ref{rmk_lost_meas}). This problem does not exist in market without transaction cost since there is no need to redefine the stock price, nor does it matter when there is a dominating measure $P$, since we can always modify a universally measurable map on a $P$-null set to make it Borel measurable. Finally, for a map $f$ on $\Omega_{t+1}$, we will often see it as a map on $\Omega_t\times \Omega_1$ and write $f=f(\omega,\omega')$. 

Define processes $X,Y$ recursively by $X_\mathbf{T}=\lbar{S}_\mathbf{T}, Y_\mathbf{T}=\ubar{S}_\mathbf{T}$ and
\begin{equation}\label{XY_recursion}
\begin{aligned}
X_t(\omega):&=\left(\inf \text{supp}_{\mathcal{P}_t(\omega)}X_{t+1}(\omega, \cdot)\right)\vee \lbar{S}_t(\omega),\\
Y_t(\omega):&=\left(\sup \text{supp}_{\mathcal{P}_t(\omega)}Y_{t+1}(\omega, \cdot)\right)\wedge \ubar{S}_t(\omega),
\end{aligned}
\end{equation}
for $t=\mathbf{T}-1, \ldots, 0$, $\omega\in\Omega_t$.

\begin{lemma}\label{XY_semi-analytic}
For each $t$, $X_t$ is l.s.a. and $Y_t$ is u.s.a..
\end{lemma}
\proof{Proof.}
We only show the lower semi-analyticity of $X_t$. A symmetric argument gives the upper semi-analyticity of $Y_t$. $X_\mathbf{T}$ is l.s.a. by assumption. Suppose $X_{t+1}$ is l.s.a., we deduce the lower semi-analyticity of $X_t$. Let $a\in\mathbb{R}$, we have
\begin{align*}
&\{\omega\in\Omega_t: \inf \text{supp}_{\mathcal{P}_t(\omega)}X_{t+1}(\omega, \cdot)<a\}\\
&=\left\{\omega\in\Omega_t: P(X_{t+1}(\omega,\cdot)<a)>0 \mbox{ for some }P\in\mathcal{P}_t(\omega)\right\}\\
&=\text{proj}_{\Omega_t}\left(\left\{(\omega,P)\in\Omega_t\times\mathfrak{P}(\Omega_1): E^P[1_{\{X_{t+1}(\omega,\cdot)<a\}}]>0\right\}\cap \text{graph}(\mathcal{P}_t)\right)
\end{align*}
It can be checked that if $A$ is an analytic set, then $1_{A}$ is an u.s.a. function. Hence $(\omega, \omega')\mapsto 1_{\{X_{t+1}(\omega, \omega')<a\}}$ is u.s.a. by the induction hypothesis. By \cite[Proposition 7.48]{StoOptCtrl}, we also have $(\omega, P)\mapsto E^P[1_{\{X_{t+1}(\omega,\cdot)<a\}}]$ is u.s.a.. It follows that $\{\omega\in\Omega_t: \inf \text{supp}_{\mathcal{P}_t(\omega)}X_{t+1}(\omega, \cdot)<a\}$ is the projection of the intersection of two analytic sets, thus also analytic. This shows $\inf \text{supp}_{\mathcal{P}_t(\omega)}X_{t+1}(\omega, \cdot)$ is l.s.a.. $X_t$, being the maximum of two l.s.a. functions, is also l.s.a..\Halmos
\endproof

\begin{remark}\label{rmk_lost_meas}
Neither Borel nor universal measurability are preserved under the backward recursion. To see this, similar to Remark 4.4 of Bouchard and Nutz \cite{BN13}, consider $\Omega_1=[0,1]$, $\mathcal{P}_t\equiv \mathfrak{P}(\Omega_1)$ and $X_{t+1}=1_{A^c}$ for some $A\subseteq\Omega_t\times\Omega_1$. Then
$\inf\text{\emph{supp}}_{\mathcal{P}_t(\omega)}X_{t+1}(\omega, \cdot)=1_{(\text{\emph{proj}}_{\Omega_t}A)^c}(\omega)$. If $A$ is Borel, then $X_{t+1}$ is Borel measurable. But $1_{(\text{\emph{proj}}_{\Omega_t}A)^c}$ is not Borel measurable in general because the projection of a Borel set may not be Borel. Similarly, if $A$ is universally measurable, then $X_{t+1}$ is universally measurable. But $1_{(\text{\emph{proj}}_{\Omega_t}A)^c}$ is not universally measurable in general because the projection of a universally measurable set may not be universally measurable. 
\end{remark}

Apart from preserving semi-analyticity, the recursively defined $[X,Y]$-market\footnote{$[X,Y]$ refers to a multi-period market with bid price $X_t$ and ask price $Y_t$ for all $t$.} has two nice properties. First, its spread is not too wide: at least all points in the interior of $[X_t, Y_t]$ admits a martingale extension to the next period $\mathcal{P}$-q.s., although there are delicate issues when the point lies on the boundary of the spread. Second, its spread is not too narrow either, in the sense that it still satisfies NA$(\mathcal{P})$ when the original market does. In summary, this new market fits our needs perfectly. The general idea of proving the nontrivial implication of multi-period FTAP is to replace the original market by the modified market $[X,Y]$, and do martingale extension in the modified market. Interior extension is not too hard in view of Proposition \ref{1pmartingale}; the challenging part is boundary extension. It turns out that boundary extension is possible if we avoid hitting boundaries as much as we can from the beginning. 

Before proving the main theorem, we need two crucial lemmas.

\begin{lemma}\label{NA_XY}
Let NA$(\mathcal{P})$ hold for the original market $[\lbar{S},\ubar{S}]$. Then NA$(\mathcal{P})$ also holds for the modified market $[X,Y]$. In particular, $X_t\leq Y_t$ $\mathcal{P}$-q.s. for all $t$.
\end{lemma}
\proof{Proof.}
We prove NA$(\mathcal{P})$ for the modified market by backward induction. Let 
\[\mathbf{M}^t:=\{[\lbar{S}_r,\ubar{S}_r]_{r=0, \ldots, t-1}, [X_{r},Y_{r}]_{r={t,\ldots, \mathbf{T}}}\}, \quad t=\mathbf{T}, \ldots, 0\] 
denotes the $(\mathbf{T}-t)$-th intermediate market obtained in the backward recursion procedure. We show NA$(\mathcal{P})$ for $\mathbf{M}^{t+1}$ implies that for $\mathbf{M}^{t}$. It suffices to show that for any self-financing portfolio process in $\mathbf{M}^t$, there exists a self-financing portfolio process in $\mathbf{M}^{t+1}$ with equal or better terminal position.

Let $\phi$ be a self-financing portfolio process in $\mathbf{M}^t$ with $\phi_0=0$ and $\phi_{\mathbf{T}+1}\geq 0$ $\mathcal{P}$-q.s.. Consider another portfolio process defined by $\eta_r:=\phi_r \ \forall r=0,\ldots, t$, 
\begin{align*}
\triangle \eta_{t+1}^1:&=1_{\{Y_t=\ubar{S}_t\}}(\triangle \phi_{t+1}^1)^+-1_{\{X_t=\lbar{S}_t\}}(\triangle \phi_{t+1}^1)^-,\\
\triangle \eta_{t+2}^1:&=\triangle \phi^1_{t+2}+1_{\{Y_t\neq\ubar{S}_t\}}(\triangle \phi_{t+1}^1)^+-1_{\{X_t\neq\lbar{S}_t\}}(\triangle \phi_{t+1}^1)^-,\\
\triangle\eta_{r+1}^1:&=\triangle \phi_{r+1}^1, \ r=t+2,\ldots, \mathbf{T},\\
\triangle \eta^0_{r+1}:&=-(\triangle \eta^1_{r+1})^+X_r+(\triangle \eta^1_{r+1})^-Y_r, \ r=t, \ldots, \mathbf{T}.
\end{align*}
 That is, we follow $\phi$ up to time $t-1$, stick to its stock position whenever the transaction at time $t$ can be carried out in $\mathbf{M}^{t+1}$, and postpone the transaction to time $t+1$ if it is not admissible in $\mathbf{M}^{t+1}$, and follow the stock position of $\phi$ again afterwards. Clearly, $\eta$ is self-financing in $\mathbf{M}^{t+1}$, and $\eta_{\mathbf{T}+1}^1=\phi_{\mathbf{T}+1}^1$. We want to show $\eta_{\mathbf{T}+1}^0\geq\phi_{\mathbf{T}+1}^0$ $\mathcal{P}$-q.s.. It suffices to show $\triangle\eta^0_{t+1}+\triangle\eta^0_{t+2}\geq \triangle\phi^0_{t+1}+\triangle\phi^0_{t+2}$. During the $(t+1)$-th and $(t+2)$-th periods, $\eta$ and $\phi$ are trading the same total number of shares, just at different times. So we only need to check that $\eta$ faces a trading price as favorable as, if not more favorable than the one faced by $\phi$. By our construction of $X_t$, when $X_t(\omega)\neq \lbar{S}_t(\omega)$, we must have $X_t(\omega)\leq X_{t+1}(\omega,\cdot)$ $\mathcal{P}_t(\omega)$-q.s.. Fubini's theorem implies $X_{t+1}\geq X_t$ $\mathcal{P}$-q.s. on $\{X_t\neq \lbar{S}_t\}$. Similarly, $Y_{t+1}\leq Y_t$ $\mathcal{P}$-q.s. on $\{Y_t\neq \ubar{S}_t\}$. Therefore, $\eta$ has price disadvantage only on a $\mathcal{P}$-polar set.\Halmos
\endproof

\begin{remark}\label{rmk:Nt}
Unlike Bouchard and Nutz \cite{BN13}, we do not attempt to show the set 
\[N_t:=\{\omega\in\Omega_t: \text{NA}(\mathcal{P}_t(\omega)) \text{ fails for the modified market}\}\]
is universally measurable and $\mathcal{P}$-polar. In fact, the measurability of $N_t$ depends on the measurability of the set-value maps $\omega \mapsto $supp$_{\mathcal{P}_t(\omega)}X_{t+1}(\omega, \cdot)$, supp$_{\mathcal{P}_t(\omega)}Y_{t+1}(\omega, \cdot)$. With $X_{t+1}$, $Y_{t+1}$ only known to be semi-analytic, the measurability of the support maps is questionable. Even if they are universally measurable, to show $N_t$ is $\mathcal{P}$-polar by a contradiction argument similar to \cite{BN13}, one needs to construct an arbitrage strategy and a measure under which a profit can be realized with positive probability. The construction of such a measure involves a measurable selection from e.g. the set-valued map $\Phi(\omega):=\{P\in\mathcal{P}_t(\omega): E^P[X_{t+1}(\omega, \cdot)-Y_t(\omega)]>0\}$. However, when $X_{t+1}$ is l.s.a. and $Y_t$ is u.s.a., we have that $\psi: (\omega, P)\mapsto E^P[X_{t+1}(\omega, \cdot)-Y_t(\omega)]$ is l.s.a. and $\{\psi>0\}$ is co-analytic. So graph$(\Phi)=\text{graph}(\mathcal{P}_t)\cap\{\psi>0\}$ fails to be analytic in general. By assuming $\mathcal{P}_t$ has a Borel graph, we can make graph($\Phi$) co-analytic. But we are not aware of a selection theorem that applies to co-analytic set, unless one is willing to assume $\Sigma_1^1$-determinancy\footnote{$\Sigma_1^1$-determinancy refers to the principle that every analytic game is determined. It cannot be proved in the standard ZFC axioms (Zermelo-Fraenkel with the Axiom of Choice).} 
(see Kechris \cite[Corollary 36.21]{Descriptive-Set-Theory}).
\end{remark}

\begin{lemma}\label{1pselection}
Let $t\in\{0, \ldots, \mathbf{T}-1\}$ and $P(\cdot):\Omega_t\rightarrow\mathfrak{P}(\Omega_1)$, $\tilde{S}_t(\cdot):\Omega_t\rightarrow \mathbb{R}$ be Borel. Let
\begin{equation}
\begin{aligned}\label{Xi}
\Xi_t(\omega):=\{&(Q_1, Q_2, \hat{P})\in\mathfrak{P}(\Omega_1)\times\mathfrak{P}(\Omega_1)\times\mathcal{P}_t(\omega): P(\omega)\ll Q_1\sim Q_2\ll \hat{P},  \\
&E^{Q_1}[Y_{t+1}(\omega, \cdot)-\tilde{S}_t(\omega)]>0 \ \text{and} \ E^{Q_2}[X_{t+1}(\omega, \cdot)-\tilde{S}_t(\omega)]<0\}, \ \omega\in\Omega_t.
\end{aligned}
\end{equation}
Then $\Xi_t$ has an analytic graph and there exist universally measurable selectors $Q_1(\cdot), Q_2(\cdot),\hat{P}(\cdot)$ for $\Xi_t$ on the universally measurable set $\{\Xi_t\neq\emptyset\}$. 
In addition, there exists a universally measurable function $\lambda(\cdot): \Omega_t\rightarrow (0,1)$ such that
\begin{equation}\label{lambda}
\lambda(\omega)E^{Q_1(\omega)}[Y_{t+1}(\omega, \cdot)]+(1-\lambda(\omega))E^{Q_2(\omega)}[X_{t+1}(\omega, \cdot)]=\tilde{S}_t(\omega) \quad \text{if} \quad \Xi_t(\omega)\neq\emptyset.
\end{equation}
\end{lemma}
\proof{Proof.}
We first show $\Xi_t$ has an analytic graph. The proof is very similar to that of Bouchard and Nutz \cite[Lemma 4.8]{BN13}. So we shall be brief. Let 
\[\Psi(\omega):=\{(Q_1,Q_2)\in\mathfrak{P}(\Omega_1)^2: E^{Q_1}[Y_{t+1}(\omega, \cdot)-\tilde{S}_t(\omega)]>0 \ \text{and} \ E^{Q_2}[X_{t+1}(\omega, \cdot)-\tilde{S}_t(\omega)]<0\}.\] 
By Bertsekas and Shreve \cite[Proposition 7.48]{StoOptCtrl}, $(\omega, Q_1, Q_2)\mapsto E^{Q_1}[Y_{t+1}(\omega, \cdot)-\tilde{S}_t(\omega)]$ is u.s.a. and $(\omega, Q_1, Q_2)\mapsto E^{Q_2}[X_{t+1}(\omega, \cdot)-\tilde{S}_t(\omega)]$ is l.s.a.. So $\Psi$ has an analytic graph. Let 
\[\Phi(\omega):=\{(Q_1, Q_2,\hat{P}) \in\mathfrak{P}(\Omega_1)^3:P(\omega)\ll Q_1\sim Q_2\ll\hat{P}\}.\] 
Define 
\[\phi(\omega,Q_1, Q_2,\hat{P}):=E^{Q_1}[dP(\omega)/dQ_1]+E^{Q_2}[dQ_1/dQ_2]+E^{Q_1}[dQ_2/dQ_1]+E^{\hat{P}}[dQ_2/d\hat{P}],\] 
where we choose a version of the Radon-Nikodym derivatives (using absolutely continuous parts) that are jointly Borel measurable (see Dellacherie and Meyer \cite[Theorem V.58]{Dellacherie82} and the remark after it). Bertsekas and Shreve \cite[Propositions 7.26, 7.29]{StoOptCtrl} then imply $\phi$ is Borel. So $\text{graph}(\Phi)=\{\phi=4\}$ is Borel. Hence $\Xi_t(\omega)=(\Psi(\omega)\times\mathcal{P}_t(\omega))\cap \Phi(\omega)$ has an analytic graph. An application of Jankov-von Neumann Selection Theorem (Lemma \ref{A-Jankov}) yields the desired universally measurable selectors $Q_1(\cdot)$, $Q_2(\cdot)$ and $\hat{P}(\cdot)$ for $\Xi_t$ on $\{\Xi_t\neq\emptyset\}$. Outside this set, define $Q_1(\cdot)=Q_2(\cdot)=\hat{P}(\cdot):=P(\cdot)$.

Next, we construct a universally measurable weight function $\lambda$. By Bertsekas and Shreve \cite[Proposition 7.46]{StoOptCtrl}, the maps 
$f:\omega\mapsto E^{Q_1(\omega)}[Y_{t+1}(\omega, \cdot)-\tilde{S}_t(\omega)]$ and $g:\omega\mapsto E^{Q_2(\omega)}[X_{t+1}(\omega, \cdot)-\tilde{S}_{t}(\omega)]$ are universally measurable. We have $f>0$ and $g<0$ on $\{\Xi\neq\emptyset\}$. Define $\lambda:=g/(g-f)$ on $\{\Xi\neq\emptyset\}$ and $\lambda:=1/2$ on $\{\Xi=\emptyset\}$. Then $\lambda$ is universally measurable, $(0,1)$-valued and satisfies \eqref{lambda}.\Halmos
\endproof

\begin{remark}\label{rmk:pair}
Lemma \ref{1pselection} is the measurable version of Proposition \ref{1pmartingale}. By Remark \ref{rmk:1p-pair},  in a one-period market, one can directly construct a CPS with the martingale component being a convex combination of the bid-ask prices. Thus it may be natural to work with
\[\Xi_t(\omega):=\{(Q,\lambda,\hat{P})\in\mathfrak{P}(\Omega_1)\times (0,1)\times\mathcal{P}_t(\omega): P(\omega)\ll Q\ll\hat{P}, E^Q[D^\lambda(\omega, \cdot)]=0\}\]
where 
$D^\lambda(\omega, \cdot):=(1-\lambda)X_{t+1}(\omega, \cdot)+\lambda Y_{t+1}(\omega, \cdot)-\tilde{S}_t(\omega)$. The problem is that when $X_{t+1}$ is l.s.a. and $Y_{t+1}$ is u.s.a., $D^\lambda$ is only universally measurable and $graph(\Xi_t)$ is not analytic in general, which makes measurable selection difficult. To overcome this issue, we break up the selection of a CPS into two steps: first select a pair of measures $(Q_1, Q_2)$ using the lower semi-analyticity of $X_{t+1}$ and the upper semi-analyticity of $Y_{t+1}$ separately, and then find a convex weight $\lambda$. A martingale can be constructed out of these selectors as an adapted convex combination of the bid-ask prices (see \eqref{ext:S}).
\end{remark}

We are now ready to prove our main results.

\proof{Proof of Theorem \ref{FTAP}}. (i) $\Rightarrow$ (ii): We first replace the original $[\lbar{S},\ubar{S}]$-market by the modified $[X,Y]$-market which lies inside $[\lbar{S},\ubar{S}]$, is still semi-analytic by Lemma \ref{XY_semi-analytic} and satisfies NA($\mathcal{P}$) by Lemma \ref{NA_XY}. It suffices to prove (ii) for the modified market because any CPS for the modified market is a CPS for the original market. Let us prove an auxiliary claim that (i) implies the following:
\begin{itemize}
\item[(ii')] $\forall P\in\mathcal{P}$, $\exists (Q,\tilde{S})\in\mathcal{Z}$, $\hat{P}\in\mathcal{P}$ with $P\ll Q\ll\hat{P}$ and an adapted process $\lambda$, valued in $[0,1]$, such that $\tilde{S}_{t}=X_t$ if $\lambda_{t-1}=0$, $\tilde{S}_{t}=Y_t$ if $\lambda_{t-1}=1$, and $\tilde{S}_{t}\in ri[X_t, Y_t]$ if $\lambda_{t-1}\in(0,1)$, $t=1, \ldots, \mathbf{T}$. Moreover, letting
\begin{align*}
\tau_1^0:&=\inf\{t\in[0,\mathbf{T}-1]:\lambda_t=0\},\\
\sigma_n^0:&=\inf\{t\in(\tau^0_n,\mathbf{T}-1]:\lambda_t>0\},\\
\tau_{n+1}^0:&=\inf\{t\in(\sigma^0_n,\mathbf{T}-1]:\lambda_t=0\},
\end{align*}
with the convention that $\inf \emptyset=\infty$, we have $X_{\tau^0_n}=Y_{\tau^0_n}=X_t$ $\forall t\in[\tau^0_n, \sigma^0_n\wedge \mathbf{T}]$ on the set $\{\tau^0_n<\infty\}$ $\mathcal{P}$-q.s.. Similarly, letting
\begin{align*}
\sigma_1^1:&=\inf\{t\in[0,\mathbf{T}-1]:\lambda_t=1\},\\
\tau^1_n:&=\inf\{t\in(\sigma^1_n,\mathbf{T}-1]:\lambda_t<1\},\\
\sigma^1_{n+1}:&=\inf\{t\in(\tau^1_n,\mathbf{T}-1]:\lambda_t=1\},
\end{align*}
we have $X_{\sigma^1_n}=Y_{\sigma^1_n}=Y_t$ $\forall t\in[\sigma^1_n, \tau^1_n\wedge \mathbf{T}]$  on the set $\{\sigma^1_n<\infty\}$ $\mathcal{P}$-q.s..
\end{itemize}

We do induction on the number of periods in the market. When there is only one period, for any $P\in\mathcal{P}$, we set $\lambda:=0$ when
$X_0=Y_0=\inf\text{supp}_{\mathcal{P}}X_1< \sup\text{supp}_{\mathcal{P}}Y_1$, $\lambda:=1$ when $X_0=Y_0=\sup\text{supp}_{\mathcal{P}}Y_1>\inf\text{supp}_{\mathcal{P}}X_1$, and $\lambda:=1/2$ when $X_0=Y_0=\inf\text{supp}_{\mathcal{P}}X_1= \sup\text{supp}_{\mathcal{P}}Y_1$. In these three cases, define $Q=\hat{P}:=P$, $\tilde{S}_0:=X_0=Y_0$ and $\tilde{S}_1:=(1-\lambda)X_1+\lambda Y_1$. It is not hard to see that NA$(\mathcal{P})$ implies $X_1=Y_0$ $\mathcal{P}$-q.s. in the first case, $Y_1=X_0$ $\mathcal{P}$-q.s. in the second case, and $X_0=Y_0=X_1=Y_1$ $\mathcal{P}$-q.s. in the third case. Hence $\tilde{S}$ is a $Q$-martingale. In all other cases (under NA$(\mathcal{P})$), we can pick $s\in (\inf\text{supp}_{\mathcal{P}}X_1, \sup\text{supp}_{\mathcal{P}}Y_1)$. Proposition \ref{1pmartingale} implies the existence of $Q_1, Q_2, \hat{P}$ such that $P\ll Q_1\sim Q_2\ll \hat{P}$ and $E^{Q_2}[X_1]<s<E^{Q_1}[Y_1]<\infty$. Exactly the same construction as in Case 3 of the proof of Proposition \ref{FTAP1periodNA} yields a desired CPS $(Q,\tilde{S})$ and a weight $\lambda\in(0,1)$. It can be checked that all the statements in (ii') are satisfied. 

Now, suppose (i) implies (ii') for any market with ${T}-1$ periods that satisfies backward recursion \eqref{XY_recursion}. We will deduce the same property for such recursively defined markets with ${T}$ periods.

Let NA($\mathcal{P}$) hold for the ${T}$-period market denoted by $\mathbf{M}$. Its submarket up to time ${T}-1$, denoted by $\mathbf{M}'$, satisfies NA($\mathcal{P}'$) where 
\[\mathcal{P}'=\{P_0\otimes\cdots\otimes P_{{T}-2}: \text{ each } P_t \text{ is a universally measurable selector of }\mathcal{P}_t\}\]
is the set of possible models for the first ${T}-1$ periods. Let $P\in\mathcal{P}$ have decomposition $P=P|_{\Omega_{{T}-1}}\otimes P_{{T}-1}$. We can apply the induction hypothesis to obtain $Q',\tilde{S}, \lambda, \hat{P}'$ described in (ii') up to time ${T}-1$ with $P|_{\Omega_{{T}-1}}\ll Q'\ll \hat{P}'\in\mathcal{P}'$. Our goal is to extend $Q',\tilde{S},\lambda,\hat{P}'$ to the ${T}$-th period.

Step 1. Extension. Lemma \ref{NA_XY} implies the universally measurable set $N:=\{X_{T-1}>Y_{T-1}\}$ is $\mathcal{P}$-polar. On $N$, we simply set $Q_{T-1}=\hat{P}_{T-1}:=P_{T-1}$, $\lambda_{T-1}:=1/2$ and $\tilde{S}_T:=(X_T+Y_T)/2$. On $N^c$, we have by the definition of $X_{T-1}, Y_{T-1}$ that
\[\inf \text{supp}_{\mathcal{P}_{T-1}}X_{T}\leq X_{T-1}\leq Y_{T-1}\leq \sup \text{supp}_{\mathcal{P}_{T-1}}Y_{T}.\] 
We perform extension on the following universally measurable sets separately:

On $A_1:=\{\inf \text{supp}_{\mathcal{P}_{T-1}}X_{T}=X_{T-1}=\tilde{S}_{T-1}<\sup \text{supp}_{\mathcal{P}_{T-1}}Y_{T}\}\cap N^c$, we set $Q_{T-1}=\hat{P}_{T-1}:=P_{T-1}$, $\lambda_{T-1}:=0$ and $\tilde{S}_T:=X_T$.

On $A_2:=\{\inf \text{supp}_{\mathcal{P}_{T-1}}X_{T}<\tilde{S}_{T-1}=Y_{T-1}=\sup \text{supp}_{\mathcal{P}_{T-1}}Y_{T}\}\cap N^c$, we set $Q_{T-1}=\hat{P}_{T-1}:=P_{T-1}$, $\lambda_{T-1}:=1$ and $\tilde{S}_T:=Y_T$.

On $A_3:=\{\inf \text{supp}_{\mathcal{P}_{T-1}}X_{T}=\tilde{S}_{T-1}=\sup \text{supp}_{\mathcal{P}_{T-1}}Y_{T}\}\cap N^c$, we set $Q_{T-1}=\hat{P}_{T-1}:=P_{T-1}$, $\lambda_{T-1}:=1/2$ and $\tilde{S}_T:=(X_T+Y_T)/2$.

On $A_4:=\{\inf \text{supp}_{\mathcal{P}_{T-1}}X_{T}<\tilde{S}_{T-1}<\sup \text{supp}_{\mathcal{P}_{T-1}}Y_{T}\}\cap N^c$, Proposition \ref{1pmartingale} guarantees the set-valued map $\Xi_{T-1}$ defined by \eqref{Xi}  is nonempty. To obtain a universally measurable selector, we first modify $P_{T-1}(\cdot)$ and $\tilde{S}_{T-1}(\cdot)$ on a $\hat{P}'$-null set $\hat{N}$ (hence $Q'$-null and $P|_{\Omega_{T-1}}$-null) to make them Borel measurable (see Bertsekas and Shreve \cite[Lemma 7.27]{StoOptCtrl}). Denote the resulting Borel kernel and random variable by $P^B_{T-1}$ and $\tilde{S}^B_{T-1}$. We can then use Lemma \ref{1pselection} to obtain universally measurable maps $Q^1_{T-1}(\cdot), Q^2_{T-1}(\cdot), \hat{P}_{T-1}(\cdot)$ and $\lambda_{T-1}(\cdot)$ such that 
$P_{T-1}^B(\omega)\ll Q^1_{T-1}(\omega)\sim Q^2_{T-1}(\omega)\ll \hat{P}_{T-1}(\omega)$, 
and if $\omega\in A_4\backslash \hat{N}$, then 
$\hat{P}_{T-1}(\omega)\in\mathcal{P}_{T-1}(\omega)$, and
\begin{equation}\label{ext:lambda}
\lambda_{T-1}(\omega)E^{Q^1_{T-1}(\omega)}[Y_{T}(\omega, \cdot)]+(1-\lambda_{T-1}(\omega))E^{Q^2_{T-1}(\omega)}[X_{T}(\omega, \cdot)]=\tilde{S}^B_{T-1}(\omega).
\end{equation}
Define
\begin{equation}\label{ext:Q}
Q_{T-1}:=\lambda_{T-1}Q^{1}_{T-1}+(1-\lambda_{T-1})Q^{2}_{T-1}.
\end{equation}
Let $Q^{1B}_{T-1}$, $Q^{2B}_{T-1}$, $Q^B_{T-1}$ be any Borel modifications of $Q^1_{T-1}$, $Q^2_{T-1}$, $Q_{T-1}$ under $\hat{P}'$. Also define
\begin{equation}\label{ext:S}
\tilde{S}_T:=\lambda_{T-1}\frac{dQ^{1B}_{T-1}}{dQ^B_{T-1}}Y_T+(1-\lambda_{T-1})\frac{dQ^{2B}_{T-1}}{dQ^B_{T-1}}X_T
\end{equation}
where we choose a version of the Radon-Nikodym derivative that are jointly Borel measurable in $(\omega, \omega')\in\Omega_t\times\Omega_1$. On $A_4\cap\hat{N}$, redefine $Q_{T-1}=\hat{P}_{T-1}:=P_{T-1}$, $\lambda_{T-1}:=1/2$ and $\tilde{S}_T:=(X_T+Y_T)/2$.

Defining $Q:=Q'\otimes Q_{T-1}$, $\hat{P}=\hat{P}'\otimes\hat{P}_{T-1}$, we have $P\ll Q\ll \hat{P}\in\mathcal{P}$ (notice that $P=P|_{\Omega_{T-1}}\otimes P_{T-1}=P|_{\Omega_{T-1}}\otimes P^B_{T-1}$), and $\tilde{S}_T\in ri[X_T, Y_T]$ $Q$-a.s..

Step 2: Verify that $E^{Q}[\tilde{S}_T|\mathcal{F}_{T-1}]=\tilde{S}_{T-1}$. This says $\tilde{S}$ is a generalized martingale, and in fact, a true martingale by Kabanov and Safarian \cite[Propositions 5.3.2, 5.3.3]{Tran.Cost}. We now check $E^{Q}[\tilde{S}_T|\mathcal{F}_{T-1}]=\tilde{S}_{T-1}$ on each $A_i$ separately.

2a) Claim that $X_T=\tilde{S}_{T-1}$ $\mathcal{P}$-q.s. (hence $Q$-a.s.) on $A_1$. Define stopping times $\tilde{\tau}_n,\tilde{\sigma}_n$ by
\begin{align*}
\tilde{\tau}_1^0:&=\inf\{t\in[0,T-2]:\lambda_t=0\},\\
\tilde{\sigma}_n^0:&=\inf\{t\in(\tilde{\tau}^0_n,T-2]:\lambda_t>0\},\\
\tilde{\tau}_{n+1}^0:&=\inf\{t\in(\tilde{\sigma}^0_n,T-2]:\lambda_t=0\}.
\end{align*}
Induction hypothesis implies that 
\begin{equation}\label{hypo}
X_{\tilde{\tau}^0_n}=Y_{\tilde{\tau}^0_n}=X_t \ \forall t\in[\tilde{\tau}^0_n, \tilde{\sigma}^0_n\wedge (T-1)]
\end{equation}
on the set $\{\tilde{\tau}^0_n<\infty\}$ $\mathcal{P}'$-q.s.. 
Define a portfolio process $\phi$ by $\phi_0:=0$, and for $t=0, \ldots, T$, 
\[\phi^1_{t+1}:=\sum_{n}1_{\{\tilde{\tau}^0_n\leq t< \tilde{\sigma}^0_n\}}-1_{\{\tilde{\tau}^0_n<\tilde{\sigma}^0_n=\infty\}\cap A_1^c\cap\{t=T-1\}}-1_{\{\tilde{\tau}^0_n<\tilde{\sigma}^0_n=\infty\}\cap A_1\cap\{t=T\}},\]
and
\[\triangle\phi^0_{t+1}=-(\triangle\phi^1_{t+1})^+X_{t}-(\triangle\phi^1_{t+1})^-Y_{t}.\]
Then $\phi$ is self-financing in the market $\mathbf{M}$. Moreover, $\phi_{T+1}\geq 0$
$\mathcal{P}$-q.s.. Indeed, on the set $\{\tilde{\tau}^0_1=\infty\}$, no trade occurs. On the set $\{\tilde{\tau}^0_m<\tilde{\sigma}^0_m=\infty\}$ for some $m\geq 1$, the strategy is to repeatedly buy a share at time $\tilde{\tau}^0_n$ and sell it at time $\tilde{\sigma}^0_n$ for all $n< m$. After that, buy a share at $\tilde{\tau}^0_m$ and sell it at time $T-1$ if $A_1$ is not observed, and at time $T$ if $A_1$ is observed. In the case where $A_1$ is not observed, \eqref{hypo} implies the selling price of every holding period is $\mathcal{P}'$-q.s. (hence also $\mathcal{P}$-q.s. since all trades occur on or before time $T-1$) the same as the buying price of that holding period. So we end up in zero position. In the case where $A_1$ is observed, \eqref{hypo} again implies $\mathcal{P}'$-q.s. perfect cancellation before the last holding period; in the last holding period, we buy a share at time $\tilde{\tau}^0_m$ at the price $Y_{\tilde{\tau}^0_m}$, and sell it at time $T$ at the price $X_T$ which is $\mathcal{P}_{T-1}$-q.s. larger than or equal to $\tilde{S}_{T-1}=X_{T-1}$ by the definition of $A_1$. \eqref{hypo} then implies $X_T\geq Y_{\tilde{\tau}^0_m}$ $\mathcal{P}$-q.s. on $\{\tilde{\tau}^0_m<\tilde{\sigma}^0_m=\infty\}\cap A_1$. So we can close our position without loss. On the set $\{\tilde{\sigma}^0_m<\tilde{\tau}^0_{m+1}=\infty\}$ for some $m\geq 1$, all trades happen on or before time $T-2$ and we have $\mathcal{P}'$-q.s. (hence $\mathcal{P}$-q.s.) perfect cancellation. By NA$(\mathcal{P})$ for the modified market (Lemma \ref{NA_XY}), we must have $\phi_{T+1}=0$ $\mathcal{P}$-q.s.. But
\begin{align*}
\phi_{T+1}=1_{A_1\cap\{\lambda_{T-2}=0\}}(X_T-\tilde{S}_{T-1}).
\end{align*}
This shows $X_T=\tilde{S}_{T-1}$ $\mathcal{P}$-q.s. on $A_1\cap\{\lambda_{T-2}=0\}$. Observe that 
\[A_1\cap\{\lambda_{T-2}>0\}\subseteq A'_1:=\{\inf \text{supp}_{\mathcal{P}_{T-1}}X_{T}=X_{T-1}=\tilde{S}_{T-1}=Y_{T-1}<\sup \text{supp}_{\mathcal{P}_{T-1}}Y_{T}\}\cap N^c.\]
We can construct another self-financing portfolio process which buys a share at time $T-1$ on the set $A'_1$ and sells it at time $T$. The resulting terminal position is $1_{A'_1}(X_T-\tilde{S}_{T-1})$ which is $\mathcal{P}$-q.s. non-negative. NA($\mathcal{P}$) then implies $X_T=\tilde{S}_{T-1}$ $\mathcal{P}$-q.s. on $A'_1$. Combining with the previous result, we have proved the claim. It follows that 
\[1_{A_1}E^Q[\tilde{S}_T|\mathcal{F}_{T-1}]=E^Q[1_{A_1}X_T|\mathcal{F}_{T-1}]=E^Q[1_{A_1}\tilde{S}_{T-1}|\mathcal{F}_{T-1}]=1_{A_1}\tilde{S}_{T-1},\]
where the first equality holds by construction, we have

2b) Similarly, we can show that $Y_T=\tilde{S}_{T-1}$ $\mathcal{P}$-q.s. on $A_2$. This implies
\[1_{A_2}E^Q[\tilde{S}_T|\mathcal{F}_{T-1}]=E^Q[1_{A_2}Y_T|\mathcal{F}_{T-1}]=E^Q[1_{A_2}\tilde{S}_{T-1}|\mathcal{F}_{T-1}]=1_{A_2}\tilde{S}_{T-1}.\]

2c) Since $X_T=\tilde{S}_{T-1}=Y_T$ $\mathcal{P}$-q.s. and $\tilde{S}_T=(X_T+Y_T)/2$ on $A_3$, we have
\[1_{A_3}E^Q[\tilde{S}_T|\mathcal{F}_{T-1}]=E^Q[1_{A_3}(X_T+Y_T)/2|\mathcal{F}_{T-1}]=E^Q[1_{A_3}\tilde{S}_{T-1}|\mathcal{F}_{T-1}]=1_{A_3}\tilde{S}_{T-1}.\]

2d) On $A_4$, using \eqref{ext:lambda}, \eqref{ext:Q} and \eqref{ext:S}, we have
\begin{align*}
E^{Q_{T-1}}[\tilde{S}_T]=E^{Q^B_{T-1}}[\tilde{S}_T]&=\lambda_{T-1}E^{Q^{1B}_{T-1}}[Y_T]+(1-\lambda_{T-1})E^{Q^{2B}_{T-1}}[X_T]\\
&=\lambda_{T-1}E^{Q^{1}_{T-1}}[Y_T]+(1-\lambda_{T-1})E^{Q^{2}_{T-1}}[X_T]\\
&=\tilde{S}^B_{T-1}=\tilde{S}_{T-1},
\end{align*}
where all equalities except the second one hold $Q'$-a.s.. This implies 
\[1_{A_4}E^Q[\tilde{S}_T|\mathcal{F}_{T-1}]=1_{A_4}\tilde{S}_{T-1}.\]

Step 3. Verify that the extended weight process $(\lambda_t)_{t=0, \ldots, T-1}$ and the corresponding stopping times satisfy the properties described in (ii'). We shall denote those stopping times for the market $\mathbf{M}$ by ${\tau}^0_n,{\sigma}^0_n,{\tau}^1_n,{\sigma}^1_n$. Notice that they differ from their counterparts $\tilde{\tau}^0_n,\tilde{\sigma}^0_n,\tilde{\tau}^1_n,\tilde{\sigma}^1_n$ for the market $\mathbf{M}'$ only possibly in the last trading cycle. We check only the properties related to ${\tau}^0_n, {\sigma}^0_n$ since our extension in step 2 is symmetric.
In this step, to keep notation simple, we treat $\mathcal{F}_t$-measurable functions as defined on $\Omega$ for each $t$, i.e. if $f$ is $\mathcal{F}_t$-measurable and $\omega|_{\Omega_t}$ is the first $t$ components of $\omega\in\Omega$, then we write $f(\omega)$ to mean $f(\omega|_{\Omega_t})$.

Let $\omega\in \{{\tau}^0_n<\infty\}$. Also assume $\omega$ belongs to the $\mathcal{P}'$-q.s. set where \eqref{hypo} hold (the $\mathcal{P}'$-q.s. set is also $\mathcal{P}$-q.s.). There are four cases.

3a) $\tau^0_n(\omega)<\sigma^0_n(\omega)\leq T-2$. In this case, $\tau^0_n(\omega)=\tilde{\tau}^0_n(\omega)$, $\sigma^0_n(\omega)=\tilde{\sigma}^0_n(\omega)$, and we have $X_{\tau^0_n(\omega)}(\omega)=Y_{\tau^0_n(\omega)}(\omega)=X_t(\omega)$ $\forall t\in[\tau^0_n(\omega),\sigma^0_n(\omega)\wedge T]=[\tilde{\tau}^0_n(\omega), \tilde{\sigma}^0_n(\omega)\wedge (T-1)]$ by the induction hypothesis.

3b) $\tau^0_n(\omega)\leq T-2$, $\sigma^0_n(\omega)= T-1$. In this case, $\tau^0_n(\omega)=\tilde{\tau}^0_n(\omega)$ and $\tilde{\sigma}^0_n(\omega)=\infty$. Again, induction hypothesis gives $X_{\tau^0_n(\omega)}(\omega)=Y_{\tau^0_n(\omega)}(\omega)=X_t(\omega)$ $\forall t\in[\tau^0_n(\omega),\sigma^0_n(\omega)\wedge T]=[\tilde{\tau}^0_n(\omega), \tilde{\sigma}^0_n(\omega)\wedge (T-1)]$.

3c) $\tau^0_n(\omega)\leq T-2$, $\sigma^0_n(\omega)=\infty$. In this case, $\tau^0_n(\omega)=\tilde{\tau}^0_n(\omega)$, $\tilde{\sigma}^0_n(\omega)=\infty$ and $[\tau^0_n(\omega),\sigma^0_n(\omega)\wedge T]=[\tilde{\tau}^0_n(\omega),  T]$. Induction hypothesis implies $X_{\tau^0_n(\omega)}(\omega)=Y_{\tau^0_n(\omega)}(\omega)=X_t(\omega)$ for $t\in[\tilde{\tau}^0_n(\omega),T-1]$. It remains to check $X_T(\omega)=X_{T-1}(\omega)$ for $\mathcal{P}$-q.s. such $\omega$. In terms of the process $\lambda$, case 3c) corresponds to $\lambda_{T-2}(\omega)=\lambda_{T-1}(\omega)=0$. $\lambda_{T-2}(\omega)=0$ implies $\tilde{S}_{T-1}(\omega)=X_{T-1}(\omega)$. Based on our construction, $\lambda_{T-1}=0$ only on $A_1$ on which we have shown that $X_T=\tilde{S}_{T-1}=X_{T-1}$ $\mathcal{P}$-q.s..

3d) $\tau^0_n(\omega)= T-1$, $\sigma^0_n(\omega)=\infty$. In this case, $\tilde{\tau}^0_n(\omega)=\tilde{\sigma}^0_n(\omega)=\infty$. We need to check $X_{T-1}(\omega)=Y_{T-1}(\omega)=X_T(\omega)$ for $\mathcal{P}$-q.s. such $\omega$. In terms of the process $\lambda$, case 3d) corresponds to $\lambda_{T-2}(\omega)>0=\lambda_{T-1}(\omega)$. By construction, $\lambda_{T-1}(\omega)=0$ only when $\omega\in A_1$. So $X_T(\omega)=\tilde{S}_{T-1}(\omega)=X_{T-1}(\omega)$ for $\mathcal{P}$-q.s. such $\omega$. If $X_{T-1}(\omega)\neq Y_{T-1}(\omega)$, then $\lambda_{T-2}(\omega)>0$ would imply $\tilde{S}_{T-1}(\omega)>X_{T-1}(\omega)$. So we must have $X_{T-1}(\omega)=Y_{T-1}(\omega)=\tilde{S}_{T-1}(\omega)=X_T(\omega)$ for $\mathcal{P}$-q.s. $\omega$ that falls into case 3d).

Statements about $\sigma^1_n, \tau^1_n$ can be verified by a symmetric argument. We therefore have proved that (i) implies (ii') for the recursively defined markets $[X,Y]$ with $T$ periods. 

Finally, we note that (ii') clearly implies (ii). 

(ii)$\Rightarrow$ (i): 
Let $f\in A_{\mathbf{T}}$ be such that $f\geq 0$ $\mathcal{P}$-q.s.. To show $f=0$ $\mathcal{P}$-q.s., we suppose on the contrary, $\exists P\in\mathcal{P}$ such that $P(\norm{f}>0)>0$ and try to derive a contradiction. Write $f=\phi_{\mathbf{T}+1}$ for some $\phi\in\mathcal{H}$. Let $(Q,\tilde{S})$ be the CPS given by (ii). According to Lemma \ref{Aequivmeas}, we can find $Q'\sim Q$ such that $\triangle \phi_{t+1}$ are $Q'$-integrable for all $t=0, \ldots, \mathbf{T}$. Then a slight modification of Kabanov and Safarian \cite[Lemma 3.2.4]{Tran.Cost} yields a bounded, strictly positive $Q'$-martingale $Z=(Z^0, Z^1)$ satisfying $\lbar{S}\leq Z^1/Z^0\leq \ubar{S}$. Let $\ab{\cdot,\cdot}$ denote the usual inner product. On one hand, 
\begin{align*}
E^{Q'}[\ab{Z_{\mathbf{T}},f}]&=\sum_{t=0}^{\mathbf{T}} E^{Q'}[\ab{Z_{\mathbf{T}},\triangle\phi_{t+1}}]=\sum_{t=0}^{\mathbf{T}} E^{Q'}[\ab{Z_t,\triangle\phi_{t+1}}]\\
&=\sum_{t=0}^{\mathbf{T}} E^{Q'}\left[Z^0_t\left(\triangle\phi^0_{t+1}+\frac{Z^1_t}{Z^0_t}\triangle\phi^1_{t+1}\right)\right]\\
&\leq \sum_{t=0}^{\mathbf{T}} E^{Q'}\left[Z^0_t\left(\triangle\phi^0_{t+1}+\ubar{S}_t(\triangle\phi^1_{t+1})^+-\lbar{S}_t(\triangle\phi^1_{t+1})^-\right)\right]\leq 0
\end{align*}
by the martingale property of $Z$ under $Q'$ and the self-financing property of $\phi$.
On the other hand, $Q'\sim Q\gg P$ implies $Q'(\norm{f}>0)>0$. Together with $f\geq 0$ and $Z_{\mathbf{T}}>0$, we get the contradictory inequality 
$E^{Q'}[\ab{Z_{\mathbf{T}},f}]>0.$\Halmos
\endproof

\section{The multi-dimensional extension}\label{sec:multi-d}
In this section, we extend our methodology and results to the multi-dimensional case. The probabilistic framework and model uncertainty set will be the same as in Section~\ref{sec:set-up}. The financial market now considered consists of $d$ dynamically traded assets where at each time $t\in\{0, \ldots, \mathbf{T}\}$, we are given a random cone $K_t:\Omega_t \twoheadrightarrow \mathbb{R}^d$, called the solvency cone. Each $K_t$ is assumed to be a closed, convex cone containing the non-negative orthant $\mathbb{R}^d_+$. It is the cone of portfolios that can be liquidated into the zero portfolio. $-K_t$ is the cone of portfolios available from zero, and $K^0_t:=K_t\cap (-K_t)$ is the linear space of portfolios that can be converted to zero and vice versa. We also define the dual cone $K^\ast_t(\omega):=\{y\in\mathbb{R}^d: \ab{x,y}\geq 0 \ \forall \ x\in K_t(\omega)\}$, $\omega\in\Omega_t$. $K^\ast_t$ is a closed, convex cone contained in $\mathbb{R}^d_+$.
We assume $graph(K^\ast_t)$ is analytic. In particular, this implies $K^\ast_t$ and $K_t$ are universally measurable by Lemmas \ref{A-analytic-graph}, \ref{A-measurable} and \ref{A-cones} (see Appendix for the definition of measurable set-valued maps).
		
Given a (random) set $G\subseteq \mathbb{R}^d$, we write $L^0(G;\mathcal{F}_t)$ for the set of all $\mathcal{F}_t$-measurable functions taking values in $G$. 
Define $A_t:=\sum_{s=0}^t L^0(-K_s; \mathcal{F}_s)$. Then $A_t$ is interpreted as the set of attainable claims at time $t$ from zero initial endowment. Also write $L^0_{\mathcal{P}}(G;\mathcal{F}_t)$ for the set of all $\mathcal{F}_t$-measurable functions taking values in $G$ $\mathcal{P}$-q.s.. 

\begin{remark}
For the market described in Section~\ref{sec:set-up}, the solvency cone $K_t$ is the closed, convex cone in $\mathbb{R}^2$ spanned by the unit vectors $e_1, e_2$ and the vectors $\ubar{S}_te_1-e_2$, $\frac{1}{\lbar{S}_t}e_2-e_1$. The self-financing condition can be equivalently written as $\triangle \phi_{t+1}\in -L^0(K_t;\mathcal{F}_t)$. 
\end{remark}

\begin{definition}
We say NA$^s(\mathcal{P})$ holds if for all $t\in\{0, \ldots, \mathbf{T}\}$,
\[A_t\cap L^0_{\mathcal{P}}(K_t; \mathcal{F}_t)\subseteq L^0_{\mathcal{P}}(K^0_t;\mathcal{F}_t).\]
\end{definition}

\begin{definition}
We say efficient friction holds if $int K^\ast_t\neq\emptyset$ $\mathcal{P}$-q.s. for all $t\in\{0, \ldots, \mathbf{T}\}$.
\end{definition}
\begin{remark}
$int K^\ast_t\neq\emptyset$ is equivalent to $K^0_t=\{0\}$. 
\end{remark}

We now compare NA$^s(\mathcal{P})$ with the no-arbitrage notion NA$_2(\mathcal{P})$ used by Bouchard and Nutz \cite{BN14}. It is completely analogous to the classical single-measure case. 
\begin{definition}[{\cite[Definition 2.2]{BN14}}]
We say NA$_2(\mathcal{P})$ holds if for all $t\in\{0, \ldots, \mathbf{T}\}$,
\[L^0_{\mathcal{P}}(K_{t+1};\mathcal{F}_t)\subseteq L^0_{\mathcal{P}}(K_{t};\mathcal{F}_t).\]
\end{definition}

\begin{lemma}\label{NA2-implies-NAs}
Assuming efficient friction, NA$_2(\mathcal{P})$ implies NA$^s(\mathcal{P})$.
\end{lemma}
\proof{Proof.}
Let $f\in A_t\cap L^0_{\mathcal{P}}(K_t;\mathcal{F}_t)$. Write $f=\xi_0+\cdots+\xi_t$ where $\xi_s\in L^0(-K_s;\mathcal{F}_s)$, $s=1, \ldots, t$. We have $\xi_0+\cdots+\xi_{t-1}\in L^0_{\mathcal{P}}(K_t;\mathcal{F}_{t-1})$ since $K_t$ is a cone. NA$_2(\mathcal{P})$ then implies $\xi_0+\cdots+\xi_{t-1}\in L^0_{\mathcal{P}}(K_{t-1};\mathcal{F}_{t-1})$. It follows that $\xi_0+\cdots+\xi_{t-2}=\xi_0+\cdots+\xi_{t-1}+(-\xi_{t-1})\in L^0_{\mathcal{P}}(K_{t-1};\mathcal{F}_{t-2})\subseteq L^0_{\mathcal{P}}(K_{t-2};\mathcal{F}_{t-2})$. Repeat this process until we have $\xi_0+\cdots+\xi_s\in L^0_{\mathcal{P}}(K_s;\mathcal{F}_s)$ for all $s=t, \ldots, 0$. Then $\xi_0\in K_0\cap (-K_0)=\{0\}$ by efficient friction. Going back one step, we have $\xi_0+\xi_1=\xi_1\in K_1\cap(-K_1)=\{0\}$ $\mathcal{P}$-q.s.. Repeat this process until we get $\xi_s=0$ $\mathcal{P}$-q.s. for all $s=0,\ldots, t$. Then $f=\xi_0+\cdots+\xi_t=0\in K^0_t$ $\mathcal{P}$-q.s..\Halmos
\endproof

\begin{remark}
In Lemma \ref{NA2-implies-NAs}, the efficient friction assumption cannot be dropped. For example, consider a simple one-period deterministic market consisting of a zero-interest rate money market and a single stock with bid-ask prices $\lbar{S}_0=\ubar{S}_0=\lbar{S}_1=1$ and $\ubar{S}_1=2$. Then $K_0=\{(x_1,x_2)\in\mathbb{R}^2: x_1+x_2\geq 0\}$ and $K_1=\{(x_1,x_2)\in\mathbb{R}^2: x_1+x_2\geq 0, x_1+2x_2\geq 0\}$. Clearly, $K_1\subseteq K_0$. So NA$_2(\mathcal{P})$ holds. But $(-1,1)\in A_1\cap L^0(K_1;\mathcal{F}_1)$ is not an element of $K^0_1=\{0\}$. So NA$^s(\mathcal{P})$ fails.
\end{remark}

\begin{remark}\label{NA2-too-strong}
The converse is not true. NA$_2(\mathcal{P})$ says a position which is not solvent today cannot be solvent tomorrow, which is quite strong. Again, consider a one-period deterministic market consisting of a zero-interest rate money market and a single stock with bid-ask prices $\lbar{S}_0=1$, $\ubar{S}_0=3$, $\lbar{S}_1=2$ and $\ubar{S}_1=4$. We have $K_0=\{(x_1,x_2)\in\mathbb{R}^2: x_1+x_2\geq 0, x_1+3x_2\geq 0\}$ and $K_1=\{(x_1,x_2)\in\mathbb{R}^2: x_1+2x_2\geq 0, x_1+4x_2\geq 0\}$. It is easy to verify by picture that this market satisfies NA$^s(\mathcal{P})$ but fails NA$_2(\mathcal{P})$. The position $(-1,2/3)$ is not solvent today, but will be solvent tomorrow, hence is an arbitrage of the second kind. We do not wish to consider this type of arbitrage because $(-1,2/3)$ is not even attainable today. 
\end{remark}

To allow for a more general market structure, we will use NA$^s(\mathcal{P})$ instead of NA$_2(\mathcal{P})$. It turns out that NA$^s(\mathcal{P})$ together with efficient friction is sufficient to imply the existence of a family of dual elements.	
	
\begin{theorem}\label{thm:FTAP-multi-d}
The following are equivalent.
\begin{itemize}
\item[(i)] NA$^s(\mathcal{P})$ and efficient friction hold.
\item[(ii)] $\forall P\in\mathcal{P}$, $\exists$ $Q\in\mathfrak{P}(\Omega)$ and a $Q$-martingale $Z$ in the filtration $(\mathcal{F}_t)$ such that $P\ll Q\lll\mathcal{P}$ and $Z_t\in int K^\ast_t$ $Q$-a.s. $\forall$ $t=0,\ldots, \mathbf{T}$.
\end{itemize}
\end{theorem}
We will refer to the above pair $(Q,Z)$ as a strictly consistent price system (SCPS).\footnote{In fact, it is a little more than an SCPS since SCPS only requires $Z_t\in ri K^\ast_t$ $Q$-a.s. for all $t$.} For the more difficult direction (i) $\Rightarrow$ (ii), our goal is to show that NA$^s(\mathcal{P})$ together with efficient friction implies that there exists a modified market with solvency cone $\widetilde{K}_t\supseteq K_t$ for all $t$, which also satisfies efficient friction and, in addition, the local version of NA$_2(\mathcal{P})$. We can then mimic the proof in Bouchard and Nutz \cite{BN14} to obtain a family of SCPS's, one for each $P\in\mathcal{P}$, for the modified market. Since $int \widetilde{K}^\ast_t\subseteq int K^\ast_t$ for all $t$, any SCPS for the modified market is also an SCPS for the original market. In the forward extension part, we cannot directly use the result of \cite{BN14} because our modified market $\widetilde{K}$ is not Borel measurable in general. We also find their assumptions of bounded friction and $K^\ast_t\cap \partial\mathbb{R}^d_+=\{0\}$ unnecessary.

\subsection{Backward recursion}

Let $\widetilde{K}^\ast_\mathbf{T}:=K^\ast_\mathbf{T}$ and for $t=\mathbf{T}-1, \ldots, 0$, $\omega\in\Omega_t$, 
\begin{equation}\label{bwd-recursion}
\widetilde{K}^\ast_t(\omega):=K^\ast_t \cap cl(conv (\Gamma_t(\omega))) \ \text{ where } \ \Gamma_t(\omega):=\text{supp}_{\mathcal{P}_t(\omega)}\widetilde{K}^\ast_{t+1}(\omega, \cdot).
\end{equation}
Here $\text{supp}_{\mathcal{P}_t(\omega)}\widetilde{K}^\ast_{t+1}(\omega, \cdot)$ denotes the smallest closed set $F\subseteq\mathbb{R}^d$ such that $P(\widetilde{K}^\ast_{t+1}(\omega, \cdot)\subseteq F)=1$ $\forall P\in\mathcal{P}_t(\omega)$. It can be shown that $y\in\text{supp}_{\mathcal{P}_t(\omega)}\widetilde{K}^\ast_{t+1}(\omega, \cdot)$ if and only if $\forall \varepsilon>0$, $\exists P\in\mathcal{P}_t(\omega)$ such that $P(\widetilde{K}^\ast_{t+1}(\omega, \cdot)\cap B_\varepsilon(y)\neq\emptyset)>0$.
Each $\widetilde{K}^\ast_t$ is well-defined, i.e. the support operation makes sense, thanks to Lemma \ref{Kstar-analytic}. $\widetilde{K}^\ast_t$ is a non-empty, closed, convex cone. Once we have $\widetilde{K}^\ast_t$, we can define the solvency cone for the modified market by $\widetilde{K}_t:=(\widetilde{K}^{\ast}_t)^\ast$. Clearly, $\widetilde{K}^\ast_t\subseteq K^\ast_t$ and $K_t\subseteq \widetilde{K}_t$. 

\begin{lemma}\label{Kstar-analytic}
For each $t$, $graph(\widetilde{K}^\ast_t)$ is analytic. 
\end{lemma}
\proof{Proof.}
By assumption, $graph(\widetilde{K}^\ast_{\mathbf{T}})=graph({K}^\ast_{\mathbf{T}})$ is analytic. Supposing $graph(\widetilde{K}^\ast_{t+1})$ is analytic, we will show $graph(\widetilde{K}^\ast_{t})$ is analytic. Observe that
\begin{align*}
graph(\Gamma_t)&=\left\{(\omega, y)\in\Omega_t\times\mathbb{R}^d: y\in \text{supp}_{\mathcal{P}_t(\omega)}\widetilde{K}^\ast_{t+1}(\omega, \cdot)\right\}\\
&=\bigcap_{n}\left\{(\omega, y)\in\Omega_t\times\mathbb{R}^d: \exists P\in\mathcal{P}_t(\omega) \text{ s.t. } P(\widetilde{K}^\ast_{t+1}(\omega, \cdot)\cap B_{1/n}(y)\neq\emptyset)>0\right\}\\
&=\bigcap_{n}\text{proj}_{\Omega_t\times \mathbb{R}^d}\Big\{(\omega, y, P)\in\Omega_t\times\mathbb{R}^d\times\mathfrak{P}(\Omega_1):  \\
& \hspace{1.15in} (\omega, P)\in graph(\mathcal{P}_t), E^P\left[1_{\{\widetilde{K}^\ast_{t+1}(\omega, \cdot)\cap B_{1/n}(y)\neq\emptyset\}}\right]>0\Big\}.
\end{align*}
Since $graph(\widetilde{K}^\ast_{t+1})$ is analytic by the induction hypothesis, we know
\begin{align*}
&\left\{(\omega, \omega', y)\in\Omega_t\times\Omega_1\times\mathbb{R}^d: \widetilde{K}^\ast_{t+1}(\omega, \omega')\cap B_{1/n}(y)\neq\emptyset\right\}\\
&=\text{proj}_{\Omega_t\times\Omega_1\times\mathbb{R}^d}\left\{(\omega, \omega', y,z)\in\Omega_t\times\Omega_1\times\mathbb{R}^d\times\mathbb{R}^d: |y-z|<1/n, (\omega, \omega', z)\in graph(\widetilde{K}^\ast_{t+1})\right\}
\end{align*}
is analytic. It follows that $(\omega, y, P)\mapsto E^P\big[1_{\{\widetilde{K}^\ast_{t+1}(\omega, \cdot)\cap B_{1/n}(y)\neq\emptyset\}}\big]$ is u.s.a.. Hence $\Gamma_t$ has an analytic graph. The rest follows from Lemma \ref{A-analytic-graph}.\Halmos
\endproof

The backward recursion can also be described using $\widetilde{K}_t$ instead of $\widetilde{K}^\ast_t$, which is more convenient for deriving no-arbitrage properties. Let
\[\Lambda_t(\omega):=\{x\in\mathbb{R}^d: x\in \widetilde{K}_{t+1}(\omega, \cdot) \ \mathcal{P}_t(\omega)\text{-q.s.}\}.\]

\begin{lemma}\label{lemma1}
$\Lambda_t= \Gamma^\ast_t.$
\end{lemma}
\proof{Proof.}
$``\subseteq"$:
Let $x\in\Lambda_t(\omega)$ and $y\in \Gamma_t(\omega)$. We want to show $\ab{x,y}\geq 0$. Let $A\subseteq \Omega_1$ be the $\mathcal{P}_t(\omega)$-q.s. set such that $x\in\widetilde{K}_{t+1}(\omega, \omega')$ $\forall \omega'\in A$. $y\in \Gamma_t(\omega)=\text{supp}_{\mathcal{P}_t(\omega)}\widetilde{K}^\ast_{t+1}(\omega, \cdot)$ implies $\forall \varepsilon>0$, $\exists P_\varepsilon\in\mathcal{P}_t(\omega)$ such that $P_\varepsilon(A_\varepsilon)>0$ where $A_\varepsilon:=\{\omega'\in\Omega_1:\widetilde{K}^\ast_{t+1}(\omega, \omega')\cap B_\varepsilon(y)\neq\emptyset\}$. For each $\varepsilon>0$, since $P_\varepsilon(A)=1$ and $P_\varepsilon(A_\varepsilon)>0$, we must have $A\cap A_\varepsilon\neq\emptyset$. Let $\omega_\varepsilon\in A\cap A_\varepsilon$. We have $x\in\widetilde{K}_{t+1}(\omega, \omega_\varepsilon)$ and $\widetilde{K}^\ast_{t+1}(\omega, \omega_\varepsilon)\cap B_\varepsilon(y)\neq\emptyset$. Pick $z_\varepsilon\in \widetilde{K}^\ast_{t+1}(\omega, \omega_\varepsilon)\cap B_\varepsilon(y)$. We have constructed a sequence $(z_\varepsilon)$ satisfying $\ab{x,z_\varepsilon}\geq 0$ and $z_\varepsilon\rightarrow y$ as $\varepsilon\rightarrow 0$. This shows $\ab{x,y}\geq 0$.

$``\supseteq"$:
For each $\omega\in\Omega_t$, consider the universally measurable set $A_\omega:=\{\omega'\in\Omega_1: \widetilde{K}^\ast_{t+1}(\omega,\omega')\subseteq \Gamma_t(\omega)\}$. By the definition of the support of set-valued maps, we have $P(A_\omega)=1$ for all $P\in\mathcal{P}_t(\omega)$. Let $x\in\Gamma^\ast_t(\omega)$. To show $x\in\Lambda_t(\omega)$, it suffices to show $x\in\widetilde{K}_{t+1}(\omega, \omega')$ $\forall \omega'\in A_\omega$. This is obvious since $\Gamma_t^\ast(\omega)\subseteq (\widetilde{K}^{\ast}_{t+1})^\ast(\omega, \omega')=\widetilde{K}_{t+1}(\omega, \omega')$ $\forall \omega'\in A_\omega$.\Halmos
\endproof

\begin{proposition}\label{prop:bwd-recursion}
\begin{itemize}
\item[(i)] $\widetilde{K}^\ast_t=K^\ast_t \cap \Lambda_t^{\ast}$.
\item[(ii)] $\widetilde{K}_t=cl(K_t+\Lambda_t).$\footnote{Given two sets $A$ and $B$ in $\mathbb{R}^d$, $A+B:=\{x+y:x\in A, y\in B\}$ denotes their Minkowski sum.}
If $int \widetilde{K}^\ast_t\neq\emptyset$, then $\widetilde{K}_t=K_t+\Lambda_t$.
\end{itemize}
\end{proposition}
\proof{Proof.}
(i) By Lemmas \ref{A-cones} and \ref{lemma1}, equation \eqref{bwd-recursion} can be rewritten as $\widetilde{K}^\ast_t=K^\ast_t \cap (\Gamma_t^{\ast})^\ast=K^\ast_t \cap \Lambda_t^{\ast}$.

(ii) Take dual on both sides of (i) and use Lemma \ref{A-cones}.\footnote{Note that the interior of finite intersection equals the intersection of interiors.}\Halmos
\endproof
\begin{corollary}\label{Ktilde-NA2}
NA$_2(\mathcal{P}_t(\omega))$ holds for the one-period market $\{\widetilde{K}_t(\omega), \widetilde{K}_{t+1}(\omega, \cdot)\}$.
\end{corollary}
\proof{Proof.}
The statement holds if and only if $\Lambda_t(\omega)\subseteq\widetilde{K}_t(\omega)$, which follows immediately from Proposition \ref{prop:bwd-recursion}(ii).\Halmos
\endproof

The next goal is to verify the modified market satisfies efficient friction. To do this, we prove efficient friction and strict no-arbitrage property at the same time using a backward induction.

\begin{lemma}\label{lemma2}
Any $f\in L^0_{\mathcal{P}}({K}_t+\Lambda_t;\mathcal{F}_t)$ admits a decomposition $f=g+h$ $\mathcal{P}$-q.s. for some $g\in L^0(K_t;\mathcal{F}_t)$ and $h\in L^0(\widetilde{K}_{t+1};\mathcal{F}_t)$.
\end{lemma}
\proof{Proof.}
Applying Filippov's implicit function theorem Himmelberg \cite[Theorem 7.1]{Himmelberg75} with $F(\omega, x,y)=x+y$ as the Carath\'{e}odory function, we get the existence of $g\in L^0(K_t;\mathcal{F}_t)$ and $h\in L^0(\Lambda_t;\mathcal{F}_t)$ such that $f=g+h$ $\mathcal{P}$-q.s.. By definition of $\Lambda_t$ and Fubini's theorem, $h\in \widetilde{K}_{t+1}$ $\mathcal{P}$-q.s..\Halmos
\endproof

\begin{proposition}\label{prop:EF}
Let NA$^s(\mathcal{P})$ and efficient friction hold for the original market. Then they also hold for the modified market.
\end{proposition}
\proof{Proof.}
Let $\mathbf{M}^t:=\{K_0, \ldots, K_{t-1}, \widetilde{K}_t, \ldots, \widetilde{K}_\mathbf{T}\}$ denote the $(\mathbf{T}-t)$-th intermediate market obtained in the backward recursion procedure. Suppose $int \widetilde{K}^\ast_r\neq\emptyset$ $\mathcal{P}$-q.s. for $r=t+1,\ldots, \mathbf{T}$, and NA$^s(\mathcal{P})$ holds for $\mathbf{M}^{t+1}$. We want to show $int \widetilde{K}^\ast_t\neq\emptyset$ $\mathcal{P}$-q.s., and NA$^s(\mathcal{P})$ holds for $\mathbf{M}^{t}$. 

Step 1. Let us show that $int \Lambda^{\ast}_t\neq\emptyset$ $\mathcal{P}$-q.s.. We have
\[\Lambda^{\ast}_t(\omega)=(\Gamma^{\ast}_t(\omega))^\ast\supseteq \Gamma_t(\omega)=\text{supp}_{\mathcal{P}_t(\omega)}\widetilde{K}^\ast_{t+1}(\omega, \cdot)\supseteq \widetilde{K}^\ast_{t+1}(\omega, \cdot) \ \mathcal{P}_t(\omega)\text{-q.s.}, \]
where the last inclusion holds by the definition of support of a set-valued map. Fubini's theorem implies $\Lambda^{\ast}_t\supseteq \widetilde{K}^\ast_{t+1}$ $\mathcal{P}$-q.s.. It then follows from the induction hypothesis that $int \Lambda^{\ast}_t\supseteq int\widetilde{K}^\ast_{t+1}\neq\emptyset$ $\mathcal{P}$-q.s..

Step 2. Now, we will verify that $int \widetilde{K}^\ast_t\neq\emptyset$ $\mathcal{P}$-q.s.. By Proposition \ref{prop:bwd-recursion}(i), it suffices to show that
\[int K^\ast_t\cap int \Lambda^{\ast}_t\neq\emptyset \quad \mathcal{P}\text{-q.s.}.\] 
Both $int K^\ast_t$ and $int \Lambda^{\ast}_t$ are $(\mathcal{P}$-q.s.) non-empty, open, convex cones. For each $\omega$ such that $int K^{\ast}_t(\omega)\neq\emptyset$ and $int \Lambda^{\ast}_t(\omega)\neq\emptyset$, but $int K^\ast_t(\omega) \cap int \Lambda^\ast_t(\omega)=\emptyset$, we can use Hahn-Banach separation theorem to obtain an $x\in\mathbb{R}^d\backslash\{0\}$ such that
$\ab{x,y}<0$ $\forall \ y\in int K^\ast_t(\omega)$ and $\ab{x,z}\geq 0$ $\forall \ z\in int \Lambda^{\ast}_t(\omega)$. The first inequality and that $int K^\ast_t(\omega)\neq\emptyset$ imply $-x\in (int K^\ast_t(\omega))^\ast=(ri K^\ast_t(\omega))^\ast=(K^{\ast}_t(\omega))^\ast=K_t(\omega)$, where we also used Lemma \ref{A-cones}. Similarly, the second inequality implies $x\in (int \Lambda^{\ast}_t(\omega))^\ast=\Lambda_t(\omega)$. 
We have therefore shown that $(-K_t(\omega)\cap \Lambda(\omega))\backslash\{0\}\neq\emptyset$.
 
Since $K^\ast_t$ and $\Gamma_t$ have analytic graphs, we know from Lemmas \ref{A-analytic-graph} and \ref{A-measurable} that $K_t=(K^\ast_t)^\ast$ and $\Lambda_t=\Gamma^\ast_t$ are universally measurable. Hence $-K_t\cap \Lambda_t$ is also universally measurable by Himmelberg \cite[Corollary 4.2]{Himmelberg75}, and has an $\mathcal{F}_t\otimes\mathcal{B}(\mathbb{R}^d)$-measurable graph by Lemma \ref{A-measurable}. It is easy to see that the graph of $(-K_t\cap \Lambda_t)\backslash\{0\}$ is also $\mathcal{F}_t\otimes\mathcal{B}(\mathbb{R}^d)$-measurable.  We can then use Lemma \ref{A-Projection} to get a universally measurable selector $x(\cdot)$ of $(-K_t\cap \Lambda_t)\backslash\{0\}$ on the set $\{(-K_t\cap \Lambda_t)\backslash\{0\}\neq\emptyset\}\supseteq \{int K^\ast_t\cap int \Lambda^\ast_t=\emptyset\}\cap\{int K^\ast_t\neq\emptyset\}\cap\{int \Lambda^\ast_t\neq\emptyset\}$. Outside this set, we define $x:=0$. Then such an $x$ belongs to $L^0(-K_t\cap\Lambda_t;\mathcal{F}_t)$. 

Now, $x\in L^0(-K_t;\mathcal{F}_t)$ implies $x\in A^{\mathbf{M}^{t+1}}_{t+1}$ where $A^{\mathbf{M}}_r$ denotes the set of attainable claims at time $r$ from zero initial endowment in a given market $\mathbf{M}$, and $x\in \Lambda_t$ implies $x\in\widetilde{K}_{t+1}$ $\mathcal{P}$-q.s. by Fubini's theorem. By the strict no-arbitrage property in the induction hypothesis, we have $x\in\widetilde{K}^0_{t+1}$ $\mathcal{P}$-q.s.. By the efficient friction property in the induction hypothesis, we have $\widetilde{K}^0_{t+1}=\{0\}$ $\mathcal{P}$-q.s.. So $x=0$ $\mathcal{P}$-q.s.. Since $x\neq 0$ on $\{int K^\ast_t\cap int\Lambda^{\ast}_t=\emptyset\}\cap\{int K^\ast_t\neq\emptyset\}\cap\{int \Lambda^{\ast}_t\neq\emptyset\}$ by our construction, and $int K^\ast_t\neq\emptyset$, $int \Lambda^{\ast}_t\neq\emptyset$ $\mathcal{P}$-q.s., $\{int K^\ast_t\cap int\Lambda^{\ast}_t=\emptyset\}$ must be $\mathcal{P}$-polar. 

Step 3. In this step, we will demonstrate that NA$^s(\mathcal{P})$ holds for $\mathbf{M}^t$.

Let $r\in\{t, \ldots \mathbf{T}\}$ and $f\in A^{\mathbf{M}^t}_r\cap L^0_{\mathcal{P}}(\widetilde{K}_r;\mathcal{F}_r)$. We want to show $f\in L^0_{\mathcal{P}}(\widetilde{K}^0_r;\mathcal{F}_r)$. 
(The case when $r\in\{0, \ldots, t-1\}$ is trivial since the solvency cones before time $t$ have not been modified in the market $\mathbf{M}^t$.) To this end, write $f=\xi_0+\cdots+\xi_r$ where $\xi_s\in L^0(-K_s;\mathcal{F}_s)$ for $s=0, \ldots, t-1$ and $\xi_s\in L^0(-\widetilde{K}_s;\mathcal{F}_s)$ for $s=t, \ldots, r$. By Step 2 and Proposition \ref{prop:bwd-recursion}, we have $\widetilde{K}_t=K_t+\Lambda_t$ $\mathcal{P}$-q.s..

Case 1. $r\geq t+1$. In this case, Lemma \ref{lemma2} implies that $\xi_t=\zeta_t+\eta_t$ $\mathcal{P}$-q.s. for some $\zeta_t\in L^0(-K_t;\mathcal{F}_t)$ and $\eta_t\in L^0(-\widetilde{K}_{t+1};\mathcal{F}_t)$. We have
\[f=\xi_0+\cdots+\xi_{t-1}+\zeta_t+(\eta_t+\xi_{t+1})+\cdots+\xi_r\in A^{\mathbf{M}^{t+1}}_r.\]
Since $f\in \widetilde{K}_r$ $\mathcal{P}$-q.s., we can use NA$^s(\mathcal{P})$ for $\mathbf{M}^{t+1}$ to obtain $f\in \widetilde{K}^0_r$ $\mathcal{P}$-q.s..

Case 2. $r=t$. In this case, first notice that $\xi_0+\cdots+\xi_{t-1}=f-\xi_t\in \widetilde{K}_t$ $\mathcal{P}$-q.s.. Lemma \ref{lemma2} implies $\xi_0+\cdots+\xi_{t-1}=g+h$ $\mathcal{P}$-q.s. for some $g\in L^0(K_t;\mathcal{F}_{t})$ and $h\in L^0(\widetilde{K}_{t+1};\mathcal{F}_{t})$. It follows that
\[h=\xi_0+\cdots+\xi_{t-1}-g\in A^{\mathbf{M}^{t+1}}_{t+1}.\]
By induction hypothesis, $h\in \widetilde{K}^0_{t+1}=\{0\}$ $\mathcal{P}$-q.s.. So $\xi_0+\cdots+\xi_{t-1}=g\in K_t$ $\mathcal{P}$-q.s.. Since
\[\xi_0+\cdots+\xi_{t-1}\in A^{\mathbf{M}^{t+1}}_{t}.\]
Induction hypothesis again yields $\xi_0+\cdots+\xi_{t-1}\in K^0_t=\{0\}$ $\mathcal{P}$-q.s.. Therefore, we have $f=\xi_t\in (-\widetilde{K}_t)\cap\widetilde{K}_t=\widetilde{K}^0_t$ $\mathcal{P}$-q.s..\Halmos
\endproof

\subsection{Forward extension}
We first state a one-period result.

\begin{lemma}\label{multi-d-1p}
Suppose $\mathbf{T}=1$ and NA$_2(\mathcal{P})$ holds. Let $P\in\mathcal{P}$. We have $\forall z\in int \widetilde{K}^\ast_0$, $\exists Q\in\mathfrak{P}(\Omega)$ and $Z_1\in L^0_{\mathcal{P}}(int \widetilde{K}^\ast_1;\mathcal{F}_1)$ such that $P\ll Q\lll\mathcal{P}$ and $E^Q[Z_1]=z$.
\end{lemma}
\proof{Proof.}
 The proof is a slight modification of that of Bouchard of Nutz \cite[Proposition 3.1]{BN14}. Let $P\in\mathcal{P}$. The first step is to show the set
$\Theta_P:=\{E^R[Y]: P\ll R\lll\mathcal{P}, Y\in L^0_{\mathcal{P}}(int\widetilde{K}^\ast_1;\mathcal{F}_1)\}$ is convex and has non-empty interior. The second step is to use separating hyperplane theorem and NA$_2(\mathcal{P})$ to show $int \widetilde{K}^\ast_0\subseteq\Theta_P$. The main change is: when showing $int\Theta_P\neq\emptyset$, due to the lack of Borel measurability of $\widetilde{K}^\ast_1$, we cannot directly select a Borel measurable $Y\in L^0_{\mathcal{P}}(int \widetilde{K}^\ast_1;\mathcal{F}_1)$, but only a universally measurable one using the analyticity of $graph(int \widetilde{K}^\ast_1)$. Nevertheless, since we are interested in showing $E^R[Y]\in int\Theta_P$ for a fixed measure $R$ satisfying $P\ll R\lll\mathcal{P}$, we can first modify $Y$ on an $R$-null set to make it Borel measurable, and then select $ \varepsilon_i \in L^0((0,1);\mathcal{F}_1), i=1, \ldots, d$ such that $Y\pm \varepsilon_i e_i\in int \widetilde{K}^\ast_1$ $R$-a.s.. Outside this $R$-a.s. set, define $\varepsilon_i:=0$. Then $E^R[\varepsilon_i]>0$, $Y\pm \varepsilon_i e_i\in L^0_{\mathcal{P}}(int \widetilde{K}^\ast_1;\mathcal{F}_1)$ and $E^R[Y\pm \varepsilon_i e_i]=E^R[Y]\pm E^R[\varepsilon_i]e_i\in \Theta_P$. We conclude that $E^R[Y]\in int (conv \{E^R[Y\pm\varepsilon_i e_i]: i=1,\ldots d\})\subseteq int\Theta_P$.\Halmos
\endproof

Let $Z$ be a one-period $Q$-martingale satisfying $Z_0=z\in int \widetilde{K}^\ast_0$ and $Z_1\in int \widetilde{K}^\ast_1$ $Q$-a.s., we can associate to $Z_1$ a vector $\mu=(\mu^1,\ldots, \mu^d)$ of probability measures on $\Omega$ defined by $d\mu^i/dQ:=Z^i_1/z^i$. Conversely, given $z\in int \widetilde{K}^\ast_0$ and $\mu=(\mu^1,\ldots, \mu^d) \in \mathfrak{P}(\Omega)^d$ satisfying $\mu\sim Q$ (i.e. $\mu^i\sim Q$ for all $i$) and $z d\mu/dQ\in int\widetilde{K}^\ast_1$ $Q$-a.s.,\footnote{Here and in the sequel, for $x,y\in\mathbb{R}^d$, $xy\in\mathbb{R}^d$ denotes their component-wise product. If we want to refer to their inner product, we will use $\ab{x,y}$.} we can define a $Q$-martingale $Z$ by $Z_0:=z$ and $Z_1:=z d\mu/dQ$. At each step in the forward extension, instead of selecting $(Q,Z_1)$, we will select $(Q, \mu)$. 

\begin{lemma}\label{lemma:multi-d-1p-selection}
Let $t\in\{0, \ldots, \mathbf{T}-1\}$ and $P(\cdot):\Omega_t\rightarrow \mathfrak{P}(\Omega_1)$, $Z_t(\cdot):\Omega_t\rightarrow\mathbb{R}^d$ be Borel. Let
\begin{equation}\label{multi-d-Xi}
\begin{aligned}
\Xi_t(\omega):=\bigg\{&(Q, \mu, \hat{P})\in \mathfrak{P}(\Omega_1)^{1+d}\times \mathcal{P}_t(\omega): P(\omega)\ll Q\sim \mu\ll \hat{P}, Z_t(\omega)\frac{d\mu}{dQ}(\cdot)\in int \widetilde{K}^\ast_{t+1}(\omega, \cdot) \ Q\text{-a.s.}\bigg\}.
\end{aligned}
\end{equation}
Then $\Xi_t$ has analytic graph and admits universally measurable selectors $Q(\cdot), \mu(\cdot), \hat{P}(\cdot)$ on the universally measurable set $\{\Xi\neq\emptyset\}$. 
\end{lemma}
\proof{Proof.}
\begin{align*}
graph(\Xi_t)=\bigg\{&(\omega, Q,\mu, \hat{P})\in\Omega_t\times \mathfrak{P}(\Omega_1)^{1+d}\times \mathcal{P}_t(\omega): P(\omega)\ll Q\sim \mu\ll \hat{P}, \\
& E^Q\left[1_{\left\{Z_t(\omega)(d\mu/dQ)(\cdot)\in int\widetilde{K}^\ast_{t+1}(\omega, \cdot)\right\}}\right]\geq 1\bigg\}
\end{align*}
For each $i$, choose a version of $\frac{d\mu^i}{dQ}(\omega')$ that is jointly Borel measurable in $(\omega', Q, \mu^i)$ by Dellacherie and Meyer \cite[Theorem V.58]{Dellacherie82}. Then the map $(\omega, \omega', Q, \mu)\mapsto Z_t(\omega)\frac{d\mu}{dQ}(\omega')$ is Borel. Since $int \widetilde{K}^\ast_{t+1}$ has an analytic graph by Lemma~\ref{A-analytic-graph}, it can be shown that the set $A:=\{(\omega, 
\omega', Q, \mu)\in \Omega_t\times\Omega_1\times \mathfrak{P}(\Omega_1)^{1+d}:Z_t(\omega)\frac{d\mu}{dQ}(\omega')\in int\widetilde{K}^\ast_{t+1}(\omega, \omega') \}$ is analytic. This implies $1_A$ and thus $(\omega, Q,\mu)\mapsto E^Q[1_A]$ are u.s.a.. Therefore, $\{(\omega, Q, \mu)\in\Omega_t\times\mathfrak{P}(\Omega_1)^{1+d}: E^Q[Z_t(\omega)\frac{d\mu}{dQ}(\cdot)\in int\widetilde{K}^\ast_{t+1}(\omega, \cdot)]\geq 1\}$ is analytic. The rest of the proof is similar to that of Lemma \ref{1pselection} or Bouchard and Nutz \cite[Lemma 4.8]{BN13}.\Halmos
\endproof

\subsection{Proof of Theorem \ref{thm:FTAP-multi-d}} \

(i)$\Rightarrow$ (ii): First construct the modified market $\widetilde{K}$ through the backward recursion \eqref{bwd-recursion}. For each $t$, $int \widetilde{K}^\ast_t\neq\emptyset$ $\mathcal{P}$-q.s. by Proposition \ref{prop:EF}. Since $\widetilde{K}^\ast_t \subseteq K^\ast_t$ for all $t$, it suffices to construct SCPS's in the modified market. Let $P\in\mathcal{P}$ be given. Pick $Z_0\in int \widetilde{K}^\ast_0\neq\emptyset$. Suppose $P=P|_{\Omega_{t-1}}\otimes P_t\otimes\cdots \otimes P_{T-1}$, and we have already constructed $Q', \hat{P}'$ up to time $t-1$ and $Z$ up to time $t$ such that $P|_{\Omega_{t-1}}\ll Q'\ll \hat{P}'$ and $Z_r\in int K^\ast_{r}$ $Q'$-a.s. for all $r\leq t$. We now extend them to the next time period in a measurable way. By modifying $P_t(\cdot)$, $Z_t(\cdot)$ on a $\hat{P}'$-null set if necessary, we may assume they are Borel measurable.

Let $A:=\{\omega\in\Omega_t: Z_t\in int \widetilde{K}^\ast_t\}\in\mathcal{F}_t$. We have $Q'(A)=1$ by construction. By Corollary \ref{Ktilde-NA2} and Lemma \ref{multi-d-1p}, we know the set-valued map $\Xi_t$ defined in \eqref{multi-d-Xi} is non-empty on $A$. Using Lemma \ref{lemma:multi-d-1p-selection}, we can find universally measurable maps $Q_{t}(\cdot), \mu_t(\cdot), \hat{P}_t(\cdot)$ such that $(Q_t, \mu_t, \hat{P}_t)\in\Xi_t$ on $A$. Modify $Q_t, \mu_t$ on a $\hat{P}'$-mull set $N$ to make them Borel. Define $Z_{t+1}:=Z_t d\mu_t/dQ_t$ where we use the Borel modification of $Q_t$ and $\mu_t$ to get a vector of jointly Borel measurable Radon-Nikodym derivatives. We have $Z_{t+1}(\omega, \cdot)\in int \widetilde{K}^\ast_{t+1}(\omega, \cdot)$ $Q_t(\omega)$-a.s. and $E^{Q_t(\omega)}[Z_{t+1}(\omega, \cdot)]=Z_t(\omega)$ $\forall \omega\in A\cap N^c$. On the $Q'$-null set $A^c\cup N$, redefine $Q_t=\hat{P}_t:=P_t$ and set $Z_{t+1}$ to be any universally measurable selector of $\widetilde{K}^\ast_{t+1}$. We have $P|_{\Omega_{t-1}}\otimes P_t\ll Q'\otimes Q_t \ll\hat{P}'\otimes \hat{P}_t$, $Z_{t+1}\in int\widetilde{K}^\ast_{t+1}$ $Q'\otimes Q_t$-a.s., and $E^{Q'\otimes Q_t}[Z_{t+1}|\mathcal{F}_t]=Z_t$.

Repeat the steps until time $\mathbf{T}$. We get measures $Q=Q'\otimes Q_t\otimes \cdots \otimes Q_{\mathbf{T}-1}$ and $\hat{P}=\hat{P}'\otimes P_t\otimes \cdots \otimes P_{\mathbf{T}-1}$ satisfying $P\ll Q\ll \hat{P}\in\mathcal{P}$, and a $Q$-generalized martingale $Z$ satisfying $Z_t\in int \widetilde{K}^\ast_t$ $Q$-a.s. for all $t$. By Kabanov and Safarian \cite[Propositions 5.3.2, 5.3.3]{Tran.Cost}, $Z$ is actually a $Q$-martingale.

(ii)$\Rightarrow$ (i): First observe that (ii) implies $int K^\ast_t\neq\emptyset$ $\mathcal{P}$-q.s. for all $t$, hence efficient friction holds. Let $f\in A_t\cap L^0_{\mathcal{P}}(K_t;\mathcal{F}_t)$. To show $f\in K^0_t=\{0\}$ $\mathcal{P}$-q.s., we suppose on the contrary $\exists P\in\mathcal{P}$ such that $P(\norm{f}>0)>0$ and try to derive a contradiction. The proof is similar to that of Theorem \ref{FTAP}. We present it here for the sake of completeness. Write $f=\sum_{r=0}^t \xi_r$ for some $\xi_r\in L^0(-K_r;\mathcal{F}_r)$. Let $(Q,Z)$ be the SCPS given by (ii). By Lemma \ref{Aequivmeas}, 
we can pick $Q'\sim Q$ such that $\xi_r$ is $Q'$-integrable for all $r=0, \ldots, t$. Then by Kabanov and Safarian \cite[Lemma 3.2.4]{Tran.Cost}, there exists a bounded $Q'$-martingale $Z'$ such that $Z'_t\in ri K^\ast_t= int K^\ast_t$ $Q'$-a.s. for all $t$. On one hand, 
\[E^{Q'}[\ab{Z'_t,f}]=\sum_{r=0}^{{t}} E^{Q'}[\ab{Z'_t,\xi_r}]=\sum_{r=0}^{{t}} E^{Q'}[\ab{Z'_r,\xi_r}]\leq 0\]
by the martingale property of $Z'$ under $Q'$ and that $Z'_r\in K^\ast_r$, $\xi_r\in -K_r$.
On the other hand, $Q'\sim Q\gg P$ implies $Q'(\norm{f}>0)>0$. Together with $f\in K_t$ and $Z'_t\in int K^\ast_t\subseteq int \mathbb{R}^d_{+}$ $Q'$-a.s., we get the contradictory inequality 
$E^{Q'}[\ab{Z'_t,f}]>0.$\Halmos

%
\begin{APPENDIX}{}
 \begin{lemma}[{\cite[Theorem VII.57]{Dellacherie82}}]\label{Aequivmeas}
Let $P\in\mathfrak{P}(\Omega)$ and $f_n$ be a sequence of ($P$-a.s. finite) random variables. There exists probability measure $R\sim P$ with bounded density with respect to $P$, such that all $f_n$ are $R$-integrable.
\end{lemma}

\begin{definition}\label{defn_meas_correspondence}
Let $(\Omega,\Sigma)$ be a measureable space and $\mathcal{X}$ be a topological space. A set-valued map $\Phi:(\Omega,\Sigma)\twoheadrightarrow \mathcal{X}$ is $\Sigma$-measurable (resp. weakly $\Sigma$-measurable) if $\Phi^{-1}(B):=\{\omega\in\Omega: \Phi(\omega)\cap B\neq\emptyset\}\in\Sigma$ for each closed (resp. open) subset $B$ of $\mathcal{X}$.
\end{definition}
When $\mathcal{X}$ is a $\sigma$-compact, separable, metrizable space (e.g. $\mathbb{R}^d$) and $\Phi$ is closed-valued, measurable and weakly measurable are equivalent (see Himmelberg \cite[Theorem 3.2(ii)]{Himmelberg75}).

\begin{lemma}\label{A-analytic-graph}
Let $\Omega$ be a Polish space and $\Phi:\Omega\twoheadrightarrow \mathbb{R}^d$ be a set-valued map with an analytic graph. The following results hold.
\begin{itemize}
\item[(a)] For any Borel set $B\subseteq\mathbb{R}^d$, $\{\omega\in\Omega: \Phi(\omega)\cap B\neq\emptyset\}$ is analytic. In particular, this implies $\Phi$ is universally measurable.
\item[(b)] $conv(\Phi)$ and $cl(\Phi)$ have analytic graphs.
\item[(c)] If $\Phi$ is convex-valued, then $int(\Phi)$ has an analytic graph.
\item[(d)] If $\Psi$ is another set-valued map with an analytic graph, then $\Phi\cap\Psi$ has an analytic graph.
\end{itemize}
\end{lemma}
\proof{Proof.}
For (a), observe that
\[\{\omega\in\Omega: \Phi(\omega)\cap B\neq\emptyset\}=\text{proj}_\Omega (graph(\Phi)\cap (\Omega\times B)).\]
For (b), observe that
\begin{align*}
graph(conv(\Phi))=\bigcup_n\text{proj}_{\Omega\times\mathbb{R}^d} \Big\{ &(\omega, y, \alpha_1, \ldots, \alpha_n, y_1, \ldots, y_n)\in\Omega\times\mathbb{R}^d\times[0,1]^d\times\mathbb{R}^{nd}: \\
& (\omega, y_i)\in graph(\Phi) \ \forall i, \sum_{i}\alpha_i=1 \text{ and } y=\sum_i\alpha_iy_i\Big\},
\end{align*}
and
\begin{align*}
graph(cl(\Phi))=\bigcap_n\text{proj}_{\Omega\times\mathbb{R}^d}\left\{(\omega, y, z)\in\Omega\times\mathbb{R}^d\times\mathbb{R}^{d}: |y-z|<1/n, (\omega, z)\in graph(\Phi)\right\}.
\end{align*}
For (c), observe that if $\Phi(\omega)$ is convex, then $y\in int(\Phi(\omega))$ if and only if $\exists n\in\mathbb{N}$ such that $y\pm e_i/n\in\Phi(\omega)$ $\forall$ $i=1, \ldots, d$. Hence,
\begin{align*}
graph(int(\Phi))=\bigcup_n\text{proj}_{\Omega\times\mathbb{R}^d}\Big\{&(\omega, y, y^\pm_1,\ldots, y^\pm_d)\in\Omega\times\mathbb{R}^d\times\mathbb{R}^{2d}: \\
& y^\pm_i=y\pm e_i/n, (\omega, y^\pm_i)\in graph(\Phi), i=1, \ldots, d\Big\}
\end{align*}
is analytic. (d) follows from $graph(\Phi\cap\Psi)=graph(\Phi)\cap graph(\Psi)$.\Halmos
\endproof

\begin{lemma}[{\cite[Lemma A.1]{BN14}}]\label{A-measurable}
Let $(\Omega, \Sigma)$ be a measurable space and $\Phi:\Omega\twoheadrightarrow {\mathbb{R}^d}$ be a non-empty, closed-valued, $\Sigma$-measurable set-valued map. The following results hold.
\begin{itemize}
\item[(a)] $\Phi^\ast$ is $\Sigma$-measurable.
\item[(b)] $graph(\Phi)$ is $\Sigma\otimes\mathcal{B}(\mathbb{R}^d)$-measurable.
\end{itemize}
\end{lemma}

\begin{lemma}[Jankov von-Neumman Selection Theorem, {\cite[Theorem 5.5.2]{A-Course-on-Borel-Sets}}]\label{A-Jankov}
Let $X,Y$ be Polish spaces and $A\subseteq X\times Y$ be an analytic set. Then there exists a universally measurable function $\phi: \text{proj}_X (A)\rightarrow Y$ such that graph$(\phi)\subseteq A$. 
\end{lemma}

\begin{lemma}[{\cite[Theorem 5.5.7]{A-Course-on-Borel-Sets}}]\label{A-Projection}
Let $(X,\Sigma)$ be a measurable space with $\Sigma$ closed under the Souslin operation, Y be a Polish space, and $A\in \Sigma\otimes\mathcal{B}(Y)$. Then $\text{proj}_X (A)\in\Sigma$, and there exists a $\Sigma$-measurable function $\phi: \text{proj}_X (A)\rightarrow Y$ such that graph$(\phi)\subseteq A$. 
\end{lemma}

\begin{lemma} \label{A-cones}
Let $K, K_1, K_2$ be non-empty cones in $\mathbb{R}^d$.
\begin{itemize}
\item[(a)] $K^\ast$ is a closed, convex cone.
\item[(b)] $(K^{\ast})^\ast=cl(conv(K))$.
\item[(c)] $(ri K)^\ast=K^\ast$.
\item[(d)] $K_1\subseteq K_2$ implies $K^\ast_1\supseteq K^\ast_2$.
\end{itemize}
If in addition, $K_1, K_2$ are convex, then
\begin{itemize}
\item[(e)] $(K_1+K_2)^\ast =K^\ast_1\cap K^\ast_2$.
\item[(f)] $(cl K_1\cap cl K_2)^\ast=cl (K^\ast_1+K^\ast_2)$, and the closure operation can be omitted if $ri K_1\cap ri K_2\neq\emptyset$.
\end{itemize}
\end{lemma}
\proof{Proof.}
(a)-(d) are standard results. (e) and (f) can be found in the book by Rockafellar \cite[Corollary 16.4.2]{Convex-Analysis}.\Halmos
\endproof
\end{APPENDIX}

%
%

\section*{Acknowledgments.}
This research is supported by the National Science Foundation under grant  DMS-0955463.


\bibliographystyle{amsplain} 
\bibliography{FTAP_bib} 


\end{document}